\documentclass{amsart}
\usepackage{amssymb, latexsym}

\let\Re=\undefined
\DeclareMathOperator{\Re}{Re}
\let\Im=\undefined
\DeclareMathOperator{\Im}{Im}

\def\bi{\begin{itemize}}
\def\bs{\begin{split}}
\def\es{\end{split}}
\def\ba{\begin{align}}
\def\bas{\begin{align*}}
\def\ea{\end{align}}
\def\eas{\end{align*}}

\def\C{{\mathbb C}} 
\def\R{{{\mathbb R}}}

\def\norm[#1][#2]{\|#1\|_{#2}}
\def\ptl[#1][#2]{\frac{\partial #1}{\partial #2}}
\def\japanese[#1]{\langle #1\rangle}

\newcommand{\eps}{{\varepsilon}}

\newcommand{\TI}{{\mathfrak{I}}}

\theoremstyle{plain}
\newtheorem{theorem}{Theorem}

\newtheorem{remark}[theorem]{Remark}
\newtheorem{proposition}[theorem]{Proposition}
\newtheorem{lemma}[theorem]{Lemma}
\newtheorem{corollary}[theorem]{Corollary}

\numberwithin{equation}{section} \numberwithin{theorem}{section}

\begin{document}

\title[A Counterexample to Dispersive Estimates]%
{A counterexample to dispersive estimates for Schr\"odinger operators in higher dimensions}
\author{M. Goldberg}
\address{Caltech}
\author{M.~Visan}
\address{University of California, Los Angeles}

\vspace{-0.3in}
\begin{abstract}
In dimension $n>3$ we show the existence of a compactly supported
potential in the differentiability class $C^\alpha$, $\alpha <
\frac{n-3}2$, for which the solutions to the linear Schr\"odinger
equation in $\R^n$,
$$
-i\partial_t u = - \Delta u + Vu , \quad u(0)=f,
$$
do not obey the usual $L^1\rightarrow L^{\infty}$ dispersive estimate. This contrasts with known results in
dimensions $n \leq 3$, where a pointwise decay condition on $V$ is generally sufficient to imply
dispersive bounds.

\end{abstract}

\maketitle

\section{Introduction}

The evolution operator for the free Schr\"odinger equation, here denoted by $e^{-it\Delta}$, is subject
to a wide variety of estimates.  Functional analysis dictates that it must be an isometry on $L^2(\R^n)$
at every fixed time $t$.  Representing $e^{-it\Delta}$ as a convolution operator with the kernel $(-4\pi
i\,t)^{-\frac{n}2}e^{-i|x|^2/(4t)}$, leads to the dispersive bound
\begin{equation}\label{nde}
\norm[e^{-it\Delta}f][\infty] \leq (4\pi |t|)^{-\frac{n}2} \norm[f][1]
\end{equation}
valid for each $t \not= 0$.  Between these two estimates one already has
most of the necessary elements to verify more subtle space-time properties
of the Schr\"odinger evolution such as global Strichartz bounds.

It is natural to ask whether a perturbed operator $e^{itH}, H =
-\Delta + V$, can satisfy (up to a constant) the same $L^1 \to
L^\infty$ estimate as the free evolution. In general, it
cannot. If $H$ has point spectrum (eigenvalues), the naive
dispersive estimate \eqref{nde} fails. Indeed, for any Schwartz
function $f$ that has nonzero inner product with an eigenfunction,
$\langle e^{itH}f,f\rangle$ does not converge to zero as $t\to
\infty$. Therefore, it is a natural endeavour to prove
\begin{equation} \label{eq:dispersive}
\norm[e^{itH}P_{ac}(H)f][\infty] \leq C|t|^{-\frac{n}2}\norm[f][1],
\end{equation}
where $P_{ac}(H)$ denotes the projection onto the absolutely continuous
spectrum\footnote{For the potentials discussed here, there is no singular
continuous spectrum by the Agmon-Kato-Kuroda
Theorem~\cite[Theorem XIII.33]{RSIV}.} of $H$.

It is known that \eqref{eq:dispersive} can fail for $t$ large in the presence of a zero-energy
eigenvalue or resonance. For more details, see \cite[Theorem 10.5]{JK},~\cite[Theorem 8.2]{Jensen},
and~\cite[\S 3]{JSS}. By assuming that zero is a regular point, that is, neither an eigenvalue nor a
resonance of $H$, one can find conditions governing the decay and regularity (but not the size, or
signature) of $V$ which are known to be sufficient to imply the dispersive bound \eqref{eq:dispersive}.
These are listed below for reference.

\begin{itemize}
\item \cite{GS1} $n=1$: \quad $(1+|x|)V \in L^1(\R)$
\item \cite{Sch} $n=2$: \quad $|V(x)| \leq C(1+|x|)^{-3-\eps}$
\item \cite{Gol} $n=3$: \quad $V \in L^{\frac32-\eps}(\R^3) \cap L^{\frac32+\eps}(\R^3)$
\item \cite{JSS} $n\geq3$: \quad $\hat{V} \in L^1$ and $(1+|x|^2)^{\gamma/2}V(x)$ is a bounded operator
   on the Sobolev space $H^\nu$  for some $\gamma>n+4$ and some $\nu>0$
\end{itemize}

For a more thorough discussion of the work on this problem, see the survey \cite{Bill}.

One might extrapolate from the results in dimensions 1, 2, and 3 that a
suitable $L^p$-type condition for potentials should be sufficient in
every dimension.  The main result of this paper, Theorem \ref{T}, shows that
this is not true:  In every dimension $n > 3$, there exist continuous and
compactly supported potentials for which the dispersive estimate
\eqref{eq:dispersive} fails.

In constructing the counterexamples, we follow the approach of \cite{GS1}
and~\cite{Gol2}. Specifically, we use Stone's formula to
construct the spectral measure from the resolvent, which in turn is studied via
a finite Born series expansion (iteration of the resolvent identity).  While
we do not explicitly separate the contributions of high and low energies,
the failure of dispersive estimates in this case should be recognized as a
high-energy phenomenon.

The three dimensional analysis of \cite{GS1} relies heavily on the simple
explicit expression of the free resolvent.  The free resolvent can be
written in terms of elementary functions in all odd dimensions; however, the
expressions become increasingly unwieldy as the dimension increases. In even
dimensions, Bessel/Hankel functions are required.
The key to avoiding this morass is the introduction of certain symbol classes,
$S^{i,j}$, which capture the essential features of the free resolvent.
In particular, in dimension $n$, one must integrate by parts approximately
$(n+1)/2$ times to obtain the appropriate power of $t$; this seems quite
impossible without such a unifying tool.

In dimensions four and higher, the Green's function is rather singular at the
origin, specifically, it is not locally square integrable.  This necessitates
carrying the Born expansion much further
than in \cite{GS1}, which adds to the complexity of our proof.

Our analysis contains certain partial positive results.  To be precise,
we show that \eqref{eq:dispersive} is attained by the tail of the Born series,
taken after a finite number (depending on the dimension) of initial terms.
The question of whether $e^{itH}P_{ac}(H)$ is dispersive then
reduces to an estimate on the initial
terms in the Born series. We
construct a potential for which the sum of these terms is bounded
below by $|t|^{-\alpha}, \alpha > \frac{n}2$, at certain times $0<t<1$. In the
limit $t\to 0$, this runs contrary to the desired bound of $|t|^{-\frac{n}2}$.
The Uniform Boundedness Principle is used
to show that the worst possible limiting behaviour can be achieved.

It should again be emphasized that the non-dispersive phenomenon takes
place over extremely short times; moreover, it is a high-energy phenomenon.
Indeed by Theorem~B.2.3 of \cite{Semigr}, for any bounded compactly supported
function $\phi$, the operator $e^{itH}\phi(H)$ maps $L^1$ into $L^\infty$ uniformly
in $t$.  This is true for very general potentials, in particular those that are bounded.

A physical interpretation is that even
high-frequency waves travelling with large velocity can be
effectively scattered by a non-smooth potential. Depending on the
geometry of the potential, the first reflection may generate an
unacceptable degree of constructive interference. For the purposes
of our counterexample, ``non-smooth'' will mean that $V$ is assumed
to possess fewer than $\frac{n-3}2$ continuous derivatives.

Compare this to the smoothness conditions in \cite{JSS}, which are
sufficient to imply a dispersive bound.  In that paper a potential is
only explicitly required
to possess derivatives of order $\nu$ for some $\nu > 0$.  Indeed,
there exist numerous examples of functions satisfying all the hypotheses of
\cite{JSS}, yet which we would consider to be non-smooth.  On the other hand,
the potentials constructed in this paper are differentiable to order
$\frac{n-3}2$ but the dispersive estimate still fails.  This suggests that
while a dispersive bound may hold for all sufficiently smooth potentials (with
rapid decay at infinity), other criteria besides the number and size of
derivatives determine what happens in the absence of such strong regularity.

The additional assumption in \cite{JSS} is that $\hat{V} \in L^1$, which is
satisfied by any potential in the Sobolev space $H^{\frac{n}2+\eps}(\R^n)$.
Determining which functions of lesser regularity also have integrable Fourier
transform is a well known difficult problem.  The counterexample constructed
here is motivated by a different and explicitly geometric consideration,
the focal pattern of reflections caused by an elliptical surface.
Strictly speaking, the reflection is caused by a highly oscillatory potential
whose level sets are ellipses. When presented in this light,
it is clear that some notion of curvature and/or convexity can
also determine whether dispersive estimates remain valid.
There is still considerable room between the currently known sufficient
conditions and the negative result presented here.  We believe this middle
ground can be explored via some combination of geometric and
Fourier analysis and that these are most likely two sides of the same coin.

%
%
%
%

\section{Notes on the free resolvent}

We introduce here a class of symbols which will be relevant in the
study of the free resolvent, simplifying  both the notation and
the analysis. For $i,j\in \mathbb{Q}$, we denote by $a_{i,j}$ a
symbol belonging to the class $S^{i,j}$, i.e., a symbol that
satisfies the following estimates
$$
\left|\frac{\partial^{k}a_{i,j}(x)}{\partial x^{k}}\right|\leq
\left \{\begin{array}{lcc} c_{k} x^{i-k} \ \ \text{if} \ \ 0<x\leq 1,\\
c_{k}x^{j-k}\ \ \text{if} \ \ x>1
\end{array} \right. \qquad \forall \ k\geq 0.
$$

The calculus of these symbols is quite straightforward: the
derivative of a symbol in $S^{i,j}$ is a symbol in $S^{i-1,j-1}$
and the product of a symbol in $S^{i,j}$ with a symbol in
$S^{i',j'}$ is a symbol in $S^{i+i',j+j'}$. In particular, the
product of a symbol in $S^{i,j}$ with $x^{\alpha}$ belongs to
$S^{i+\alpha,j+\alpha}$.

Now let us consider the resolvent of the free Schr\"odinger
equation,
$$
R_{0}(z)=(-\Delta-z)^{-1}.
$$
In dimension $n\geq4$, $R_{0}(z)$ is given by the kernel:
\begin{align}\label{free kernel}
R_{0}(z)(x,y)=\frac{i}{4}\Bigl(\frac{z^{\frac{1}{2}}}{2\pi|x-y|}\Bigr)^{\frac{n}{2}-1}H_{\frac{n}{2}-1}^{(1)}(z^{\frac{1}{2}}|x-y|)
,
\end{align}
where $\Im z^{\frac{1}{2}} \geq 0$ and $H_{\frac{n}{2}-1}^{(1)}$ is the
first Hankel function.

We encode the information contained in the asymptotic expansions
of the first Hankel function near the origin and at infinity (see
\cite{GR}), together with the information provided
by the differential equation satisfied by the first Hankel
function,
$$
H_{\nu-1}^{(1)}(z) - H_{\nu+1}^{(1)}(z)=2\frac{d}{dz}
H_{\nu}^{(1)}(z),
$$
into the following formula valid for $\Re \nu>-\frac{1}{2}$ and
$|\arg z|<\pi$,
$$
H_{\nu}^{(1)}(z)=e^{iz} a_{-\nu,-\frac{1}{2}}(z).
$$
This together with \eqref{free kernel} yield a representation for
the kernel of the free resolvent in dimension $n\geq 4$ in terms
of the aforementioned symbols, that is,
\begin{align}\label{R0}
R_0^{\pm}(\lambda^2)(x,y)=a_{0,\frac{n-3}{2}}(\lambda|x-y|)\frac{e^{\pm
i\lambda|x-y|}}{|x-y|^{n-2}},
\end{align}
where $R_0^{\pm}(\lambda^2)$ denote the boundary values
$R_0(\lambda^2\pm i0)$.

Let us also point out a similar formula for the imaginary part of
the free resolvent,
\begin{align}\label{ImR0}
\Im
R_{0}(\lambda^{2})(x,y)=a_{n-2,\frac{n-3}{2}}(\lambda|x-y|)\frac{e^{\pm
i\lambda|x-y|}}{|x-y|^{n-2}},
\end{align}
by which we mean that we can write it as the sum of two terms of this type, one with phase $e^{i\lambda|x-y|}$
and the other with phase $e^{-i\lambda|x-y|}$. Indeed, using (for example) the identity
\begin{align}\label{scaling}
\lambda^{n-2}(-\Delta-1)^{-1}(\lambda x)=(-\Delta-\lambda^{2})^{-1}(x)
\end{align}
for the kernels of the free resolvents, we can write
$$
\Im
R_{0}(\lambda^{2})(x,y)=\lambda^{n-2}(\lambda|x-y|)^{\frac{2-n}{2}}J_{\frac{n-2}{2}}(\lambda|x-y|),
$$
where $J_{\frac{n-2}{2}}$ denotes the Bessel function. Consulting
the asymptotic expansions of the Bessel function near the origin
and at infinity (see again \cite{GR}) and using the
differential equation satisfied by the Bessel function,
$$
J_{\nu-1}(z) - J_{\nu+1}(z)=2\frac{d}{dz} J_{\nu}(z),
$$
one easily derives \eqref{ImR0}.

The purpose of understanding the free resolvent is that it enables us to study functions of $H$
through to the Stone formula for the spectral measure:
\begin{align*}
\bigl\langle F(H)P_{ac} f,g\bigr\rangle
    &=2\int_{0}^{\infty} F(\lambda^2)\lambda
\langle E'(\lambda^2)f, g\rangle d\lambda  \\
    &=\frac{2}{\pi} \int_{0}^{\infty}F(\lambda^2)\lambda\bigl\langle
\Im R_{V}(\lambda^2)f, g\rangle d\lambda,
\end{align*}
where $f,g$ are any two Schwartz functions, $P_{ac}$ denotes the projection onto the absolutely
continuous spectrum of $H$, $E'(\lambda)$ denotes the spectral measure associated to $H$, and
$R_V^\pm(\lambda^2) := (H-\lambda^2\pm i0)^{-1}$ is the resolvent of the perturbed Schr\"odinger
equation.  We have chosen signs so that $2i \Im R_V(\lambda^2) = R_V^+(\lambda^2) - R_V^-(\lambda^2)$.

In order to compute the kernel of $\Im R_V(\lambda^2)$, we make
use of the resolvent identity:
$$
R_{V}^{\pm}(\lambda^{2})=R_{0}^{\pm}(\lambda^{2})-R_{0}^{\pm}(\lambda^{2})
VR_{V}^{\pm}(\lambda^{2}),
$$
which by iteration gives rise to the following finite Born series
expansion:
\begin{align}
R_{V}^{\pm}(\lambda^{2})& \label{Born series} \\
   &=\sum_{l=0}^{2m+1}R_{0}^{\pm}(\lambda^{2})[-VR_{0}^{\pm}(\lambda^{2})]^{l}
 \label{first terms} \\
 &\quad +R_{0}^{\pm}(\lambda^{2})V[R_{0}^{\pm}(\lambda^{2})V]^{m}R_{V}^{\pm}
(\lambda^{2})[VR_{0}^{\pm}(\lambda^{2})]^{m}VR_{0}^{\pm}(\lambda^{2}).
\label{tail}
\end{align}
Elementary algebra can also be used to solve for $R_V^\pm(\lambda^2)$ in terms of $R_0^\pm(\lambda^2)$:
$$R_V^\pm(\lambda^2) = \big(I + R_0^\pm(\lambda^2)V\big)^{-1}R_0^\pm(\lambda^2)
  := S^\pm(\lambda^2)R_0^\pm(\lambda^2)$$
For now this identity is only a formal statement, as we have not shown that
$S^\pm(\lambda^2)=(I + R_0^\pm(\lambda^2)V\big)^{-1}$ exists as
a bounded operator on any space.  Existence and uniform boundedness of $S^\pm(\lambda^2)$ will be
demonstrated in Section~\ref{sec:tail}.

%
%
%
%

\section{Useful lemmas} \label{sec:lemmas}

In this section we prove a few technical lemmas. We begin with
certain results related to the boundedness of the Riesz potentials
between various weighted spaces. By Riesz potentials, we mean the
operators
$$
I_{\alpha}:f\mapsto|x|^{\alpha-n}*f
$$
where $0<\alpha< n$.

Let $\TI_{q}$ denote the space of compact operators $T$ for which
$\|T\|_{\TI_{q}}=[tr(|T|^{q})]^{\frac{1}{q}}$ is finite. We recall
the following well-known result (see \cite[Theorem XI.20]{RSIII}):

\begin{lemma}\label{opbounds}
Let $f,g\in L^{q}(\mathbb{R}^{n})$, for some $2\leq q<\infty$.
Then, $f(x)g(-i\nabla)\in\TI_{q}$ and
$$
\|f(x)g(-i\nabla)\|_{\TI_{q}}\leq(2\pi)^{-\frac{n}{q}}\|f\|_{q}\|g\|_{q}.
$$
Here, $f(x)$ denotes multiplication by $f$ in physical space,
while $g(-i\nabla)$ denotes multiplication by g in frequency
space.
\end{lemma}

As a consequence of Lemma \ref{opbounds}, one can derive results on the boundedness of the Riesz
potentials between various weighted spaces.  To describe these spaces, we will use the notation
$$\norm[f][L^{p,\sigma}] := \norm[\langle x\rangle^\sigma f][L^p],$$
where $\langle x\rangle :=(1+|x|^2)^{1/2}$, $1\leq p\leq \infty$, and $\sigma\in \R$. Following the notation of
Jensen and Kato, we write $B(0,\sigma;0,-\sigma')$ for the set of bounded operators from $L^{2,\sigma}$ to
$L^{2, -\sigma'}$, while $B_{0}(0,\sigma;0,-\sigma')$ denotes the set of compact operators from
$L^{2,\sigma}$ to $L^{2, -\sigma'}$, Jensen shows (see Lemma 2.3 in \cite{Jensen}) the following result.

\begin{proposition}\label{impprop}
{\rm1)} If $0<\alpha<\frac{n}{2}$, $\sigma, \sigma'\geq 0$, and $\sigma+\sigma'\geq \alpha$, then
$I_{\alpha}\in B(0,\sigma;0,-\sigma').$ Moreover, if $\sigma+\sigma'>\alpha,$ then $I_{\alpha}\in
B_{0}(0,\sigma;0,-\sigma').$

{\rm2)} If $\frac{n}{2}\leq\alpha<n$, $\sigma,\sigma'>\alpha-\frac{n}{2}$, and
$\sigma+\sigma'\geq\alpha$, then $I_{\alpha}\in B(0,\sigma;0,-\sigma').$ Moreover, if
$\sigma+\sigma'>\alpha$, then $I_{\alpha}\in B_{0}(0,\sigma;0,-\sigma').$
\end{proposition}

The case $\alpha \geq n$ may appear qualitatively different from the Riesz potentials considered above;
however, the mapping bounds between weighted $L^2$ spaces are still valid.

\begin{proposition}\label{alpha>n}
Let $\alpha \geq n$.  The convolution operator $I_\alpha := f \mapsto |x|^{\alpha-n} * f$ is an element
of $B_0(0,\sigma;0,-\sigma')$, provided $\sigma, \sigma' > \alpha - \frac{n}2$.
\end{proposition}
\begin{proof}
As every Hilbert-Schmidt operator is compact, in order to prove the proposition it suffices to show that
$I_\alpha$ is a Hilbert-Schmidt operator between $L^{2,\sigma}$ and $L^{2,-\sigma'}$. In turn, this is
equivalent to showing the finiteness of the integral
$$\iint \langle x\rangle^{-2\sigma} |x-y|^{2(\alpha-n)} \langle y\rangle^{-2\sigma'}\,dx dy.$$
Consider the integral with respect to $x$, namely
$$\int \langle x\rangle^{-2\sigma}|x-y|^{2(\alpha-n)}dx.$$

If $|y| \leq 1$, this is dominated by the integral of $\langle x\rangle^{2(\alpha-\sigma-n)}$, which is
finite because $\sigma > \alpha-\frac{n}2$.

Now suppose $|y| > 1$.  Over the region where $|x| \leq \frac12|y|$, the factor $|x-y|$ is essentially
of size $|y|$, as can be seen from the triangle inequality.  Meanwhile, the factor $\langle x\rangle^{-2\sigma}$
is integrable because $\sigma > \alpha-\frac{n}2 \geq \frac{n}2$.  Consequently, the integral over this region
is bounded by $|y|^{2(\alpha-n)}$, i.e.,
$$
\int_{|x| \leq \frac12|y|}\langle x\rangle^{-2\sigma}|x-y|^{2(\alpha-n)}dx \lesssim |y|^{2(\alpha-n)}.
$$
Over the region where $|x-y| \leq \frac12|y|$, the triangle inequality dictates $|x|\sim|y|$. Hence,
\begin{align*}
\int_{|x-y| \leq \frac12|y|}\langle x\rangle^{-2\sigma}|x-y|^{2(\alpha-n)}dx
&\lesssim \langle y \rangle ^{-2\sigma} \int_{|x-y| \leq \frac12|y|}|x-y|^{2(\alpha-n)}dx\\
&\lesssim \langle y\rangle^{-2\sigma}|y|^{2\alpha-n}
\lesssim |y|^{2\alpha-2\sigma-n}.
\end{align*}
Everywhere else in $\R^n$, the two functions $|x|$ and $|x-y|$
are of comparable size.  Recalling that $\sigma>\alpha-\frac{n}{2}$, the integral over this region is then dominated
by
$$
\int_{|x|>\frac12|y|} \langle x\rangle^{-2\sigma}|x|^{2(\alpha-n)}
\lesssim \int_{|x|>\frac12|y|}|x|^{2(\alpha-\sigma-n)}
\lesssim |y|^{2\alpha-2\sigma-n}.
$$
Therefore, the dominant term for large $y$ comes from the region $|x|\leq \frac{1}{2}|y|$.
To complete the estimate for the Hilbert-Schmidt norm, it remains to bound
the integral over the $y$-variable. As $\sigma'>\alpha-\frac{n}{2}$, this is dominated by
$$\int \langle y\rangle^{2(\alpha-n)} \langle y\rangle^{-2\sigma'}\,dy \lesssim 1.$$
This concludes the proof of Proposition~\ref{alpha>n}.
\end{proof}

Propositions \ref{impprop} and \ref{alpha>n} immediately yield some mapping bounds for
the free resolvent and its derivatives. Indeed, we have

\begin{corollary}\label{deriv in weighted}
Let $j$ be any nonnegative integer and suppose $\sigma, \sigma' > j+\frac12$
with $\sigma + \sigma' > j + \frac{n+1}2$.  Then
\begin{equation*}
\Big\| \Big(\frac{d}{d\lambda}\Big)^{j}R_0^\pm(\lambda^2)f\Big\|_{L^{2,-\sigma'}}
\lesssim \lambda^{-j}\japanese[\lambda]^{j+\frac{n-3}2} \norm[f][L^{2,\sigma}].
\end{equation*}
\end{corollary}
\begin{proof}
Recall that the kernel of $R_0^\pm(\lambda^2)$ is given by
$|x|^{2-n}e^{\pm i\lambda|x|}a_{0,\frac{n-3}2}(\lambda |x|)$.
When a symbol is differentiated, the effect is comparable to
dividing by $\lambda$; see Section~2 for the calculus of the symbols $a_{i,j}$.
Each derivative that falls on the exponential factor increases the power of $|x|$ by one.

Based on these possible outcomes, the integral kernel of
$(\frac{d}{d\lambda})^jR_0^\pm(\lambda^2)$ must be of the
form $\lambda^{-j}|x|^{2-n}e^{\pm i\lambda|x|}a_{0,\frac{n-3}2+j}(\lambda|x|)$,
which is dominated pointwise by the kernel of $\lambda^{-j}I_2 + \lambda^{\frac{n-3}2}
I_{\frac{n+1}2+j}$. Thus, for the kernels, we have the pointwise inequality
\begin{align}\label{deriv free res bound}
\Big(\frac{d}{d\lambda}\Big)^jR_0^\pm(\lambda^2)
\lesssim \lambda^{-j}\japanese[\lambda]^{j+\frac{n-3}2}(I_2 +I_{\frac{n+1}2+j}).
\end{align}
The claim follows from Propositions \ref{impprop} and \ref{alpha>n}.
\end{proof}

The estimate above is based entirely on the size of the integral kernel of
$R_0^\pm(\lambda^2)$ and its derivatives and completely ignores the
oscillatory nature of these functions.  If one takes advantage of this
oscillation using Fourier analysis techniques, the result is a much more
subtle mapping estimate known as the Limiting Absorption Principle for the
free resolvent (see \cite{agmon}, \cite[Theorem XIII.33]{RSIV}).

\begin{lemma}\label{lap}
Choose any $\sigma, \sigma' > \frac12$ and $\eps > 0$.  Then for all
$\lambda \geq 1$,
\begin{equation*}
\norm[R_0^\pm(\lambda^2)f][L^{2,-\sigma'}] \lesssim \lambda^{-1+\eps}
\norm[f][L^{2,\sigma}].
\end{equation*}
\end{lemma}
\begin{proof}[Sketch of Proof]
First, one shows that $R_0^\pm(1)$ is a bounded operator from $L^{2,\sigma}$
to $L^{2,-\sigma'}$.  One characterization of $R_0^\pm(1)$ is that it
multiplies the Fourier transform of $f$ by the distribution
$m(\xi) = \frac{c_n}{|\xi|^2-1} \pm C_n i \delta_0(|\xi|^2 - 1)$.

If $f \in L^{2,\sigma}$, then $\hat{f} \in H^\sigma(\R^n)$.  As
$\sigma > \frac12$, the Trace Theorem (see \cite{agmon} or \cite[Theorem IX.39]{RSII}) implies that
$\hat{f}$ will restrict to an $L^2$ function on surfaces of codimension 1.
The surface of particular interest here is the unit sphere, where $m(\xi)$ becomes singular.
After a partition of unity decomposition and smooth changes of variables, each sector of the sphere
can be mapped to a subset of the hyperplane $\{\xi_1 = 0\}$.  Under the
same change of variables, the singular part of $m(\xi)$ takes the form
$m(\xi) = \frac1{\xi_1} \pm i \delta_0(\xi_1)$.

This reduces matters to a one-dimensional problem.  In $\R$, multiplying
the Fourier transform by $\frac{1}{\xi_1}$ or by a delta-function are
integration operators which map $L^1$ to $L^\infty$ and consequently also
map $L^{2,\sigma}$ to $L^{2,-\sigma'}$ provided $\sigma, \sigma' > \frac12$.

The kernel of $R_0^\pm(\lambda^2)$ is simply a dilation of $R_0^\pm(1)$; see \eqref{scaling}.
A straightforward scaling argument shows that
\begin{align*}
\norm[R_0^\pm(\lambda^2)f][L^{2,-\sigma'}] \lesssim \lambda^{\sigma+\sigma'-2} \norm[f][L^{2,\sigma}]
\end{align*}
for all $\lambda \ge 1$. Finally, one can use the embeddings
$L^{2,\sigma} \subset L^{2, \min(\sigma, \frac{1+\eps}2)}$ and
$L^{2,-\min(\sigma',\frac{1+\eps}2)} \subset L^{2, -\sigma'}$
to obtain the desired power of decay in $\lambda$.
\end{proof}

Note that Corollary \ref{deriv in weighted} and Lemma \ref{lap} imply that the free resolvent and its derivatives
map functions with good decay at infinity to functions with less decay.
If this is composed with multiplication by a potential $V(x)$ with
sufficient decay at infinity, the resulting operator will be bounded
from certain weighted spaces to themselves.

\begin{corollary}
Let $j$ be a nonnegative integer and suppose
$|V(x)| \leq C\japanese[x]^{-\beta}$ for some $\beta > \max(\frac{n+1}2 + j,
2j+1)$. Then for every $j+\frac12 < \sigma < \beta-(j+\frac12)$,
\begin{equation} \label{freeest}
\Big\|\Big(\frac{d}{d\lambda}\Big)^jR_0^\pm(\lambda^2)Vf\Big\|_{L^{2,-\sigma}}
\leq \left\{ \begin{aligned} &\japanese[\lambda]^{-1+\eps}\norm[f][L^{2,-\sigma}], &&{\rm if}\ j=0, \\
&\lambda^{-j}\japanese[\lambda]^{j+\frac{n-3}2}\norm[f][L^{2,-\sigma}],
&&{\rm if}\ j \geq 1.  \end{aligned}\right.
\end{equation}
\end{corollary}
\begin{remark}
It is possible to mimic the proof of the Limiting Absorption Principle to
prove stronger estimates in the cases where $1 \leq j < \frac{n-1}2$.
These are interesting in their own right, but will not be needed here.
\end{remark}

As mentioned in the introduction, the kernel of the free resolvent is not
locally square integrable, which places it outside the context of the mapping estimates above.
However, as the next results demonstrate, the kernel
associated to $[VR_0^\pm ]^m$ belongs to a weighted $L^2$ space,
provided $m$ is big enough and $V$ decays sufficiently rapidly.

We start with the following
\begin{lemma}\label{easylem}
Let $\mu$ and $\sigma$ be such that $\mu<n$ and
$n<\sigma+\mu$. Then
$$
\int_{\mathbb{R}^{n}}\frac{dy}{\langle y
\rangle^{\sigma}|x-y|^{\mu}}\lesssim
\begin{cases}
\langle x \rangle^{n-\sigma-\mu}, & \sigma<n\\
\langle x \rangle^{-\mu}, & \sigma>n.
\end{cases}
$$
\end{lemma}

\begin{proof}

We analyze the integral on each of the following three disjoint
domains:

Domain 1: $|y| \leq \frac{|x|}{2}$. From the triangle inequality
we get $|x-y|\sim|x|$; we estimate the contribution of this domain
to the integral by
$$
|x|^{-\mu} \int_{|y|\leq \frac{|x|}{2}}\langle y\rangle ^{-\sigma}
dy\lesssim |x|^{-\mu}|x|^n \langle x\rangle^{-\sigma}
     \lesssim
\begin{cases}
\langle x \rangle^{n-\sigma-\mu}, & \sigma<n\\
\langle x \rangle^{-\mu}, & \sigma>n.
\end{cases}
$$

Domain 2:  $|x-y|\leq\frac{|x|}{2}$. On this domain $|y|\sim|x|$
and we estimate its contribution to the integral by
$$
\int_{|x-y|\leq\frac{|x|}{2}}\frac{dy}{\langle
x\rangle^{\sigma}|x-y|^{\mu}}
    =\langle x\rangle^{-\sigma} \int_{0}^{\frac{|x|}{2}}\frac{r^{n-1}}{r^{\mu}}dr
    \lesssim \langle x\rangle^{n-\sigma-\mu},
$$
where the inequality holds because $\mu<n$.

Domain 3:  $|y|>\frac{|x|}{2}$ and $|x-y|>\frac{|x|}{2}$. The
triangle inequality yields $|x-y|\sim|y|$ and as $n-\sigma-\mu<0$,
we obtain the estimate
$$
\int_{|y|>\frac{|x|}{2}}\langle y\rangle^{-\sigma}|y|^{-\mu}dy
   \lesssim \langle x\rangle^{n-\sigma-\mu},
$$
by treating $|x|\leq 1$ and $|x|>1$ separately.
\end{proof}

\begin{proposition} \label{smoothing} Suppose $|V(x)| \leq C\japanese[x]^{-\beta}$ for some
$\beta > n+3$.  Then for any integer $0 \leq j \leq \frac{n}2+1$ and any pair $(p,q)$ such that
either $1< p<\frac{2n}{n+3}$ and $\frac{1}{q}=\frac{1}{p}-\frac{2}{n}$, or $p=1$ and $1\leq q<\frac{n}{n-2}$, we have
\begin{equation*}
\Big\| V \Big(\frac{d}{d\lambda}\Big)^jR_0^\pm(\lambda^2)f\Big\|_{L^{1,\frac32}\cap L^q}
\lesssim \lambda^{-j}\japanese[\lambda]^{j+\frac{n-3}2} \norm[f][L^{1,\frac32}\cap L^{p}].
\end{equation*}
\end{proposition}
\begin{proof}
In view of \eqref{deriv free res bound}, we need only prove estimates for the operator $V I_k$
for certain $2 \leq k \leq n+\frac32$.

The weighted $L^1$ estimate follows from
$$
\sup_{x\in\R^n} \japanese[x]^{-\frac32} \int_{\R^n} \japanese[y]^{\frac32-\beta}|x-y|^{k-n}\, dy  \lesssim 1,
$$
which is a direct consequence of Lemma~\ref{easylem} with $\sigma=\beta-\frac32$ and $\mu=n-k$.

We turn now to the smoothing estimate.  Consider first the case $p=1$.
Lemma~\ref{easylem} with $\sigma=q\beta$ and $\mu=q(n-k)$ implies that for $1\leq q<\frac{n}{n-2}$, we have
\begin{align*}
\int |V(x)|^q|x-y|^{q(k-n)}dx \lesssim \langle y \rangle^{\frac{3q}{2}}.
\end{align*}
Note that the upper bound on $q$ is dictated by $k=2$.
Thus, in the case $p=1$, the claim follows from Minkowski's inequality:
\begin{align*}
\Big\| V(x)\int|x-y|^{k-n}|f(y)|dy\Big\|_{L_x^q}
&\lesssim \int|f(y)|\langle y \rangle^{\frac{3}{2}} dy
\lesssim \norm[f][L^{1,\frac32}].
\end{align*}
Lastly, we treat the case $1<p<\frac{2n}{n+3}$. Note that given $p$, the choice of $q$ is governed by the
Hardy-Littlewood-Sobolev inequality for $I_2$. As $V\in L^\infty $, we obtain
\begin{align*}
\|V I_2(f)\|_{L^q}\lesssim \|f\|_{L^p}\lesssim  \norm[f][L^{1,\frac32}\cap L^{p}].
\end{align*}
It remains to consider $I_k$ with $k=\frac{n+1}2 + j$.
For $0\leq j<\frac{n-1}{2}$, by the Hardy-Littlewood-Sobolev inequality and the fact that $V\in L^1\cap L^\infty$, we get
\begin{align*}
\|V I_{\frac{n+1}{2}+j}(f)\|_{L^q}
&\lesssim \|V\|_{\frac{2pn}{n-3+2j}}\|  I_{\frac{n+1}{2}+j}(f)\|_ {\frac{2pn}{n-1-2j}}
\lesssim  \|f\|_{L^p}\lesssim \norm[f][L^{1,\frac32}\cap L^{p}].
\end{align*}
For the remaining values of $j$, i.e., $\frac{n-1}{2}\leq j\leq \frac{n}{2}+1$, we use again
Lemma~\ref{easylem} with $\sigma=q\beta$ and $\mu=q(n-k)$ to obtain
\begin{align*}
\int |V(x)|^q|x-y|^{q(k-n)}dx \lesssim \langle y \rangle^{\frac{3q}{2}}
\end{align*}
for $1\leq q< \frac{2n}{n-1}$.  For the values of $p$ currently under consideration, $q$ is guaranteed to lie
in this range. Another application of Minkowski's inequality yields
\begin{align*}
\|V I_{\frac{n+1}{2}+j}(f)\|_{L^q}
\lesssim  \norm[f][L^{1,\frac32}] \lesssim \norm[f][L^{1,\frac32}\cap L^{p}].
\end{align*}
This completes the proof of the proposition.
\end{proof}

\begin{proposition}\label{weighted}
For any $0 \leq j \leq \frac{n}2+1$ and $\sigma > j+\frac12$,
$$ \Big\|\Big(\frac{d}{d\lambda}\Big)^j R_0^\pm(\lambda^2)f
    \Big\|_{L^{2,-\sigma}} \lesssim   \lambda^{-j}
\japanese[\lambda]^{j+\frac{n-3}2}\norm[f][L^{1,\frac32}\cap L^2]. $$
\end{proposition}
\begin{proof}
We use the estimate \eqref{deriv free res bound} and split the resolvent kernel into two pieces, according to whether
$|x-y| < 1$ or $|x-y| \geq 1$.  The piece supported away from the diagonal $x=y$ maps $L^1$
into $L^{2,-\sigma}$ because of the bound
$$ \sup_{x\in\R^n} \int_{|x-y|\geq 1}\frac{dy}{|x-y|^{2(n-k)}\langle y\rangle^{2\sigma}}
\lesssim 1, $$
valid for any $k\leq \frac{n+1}2+j$.
The piece supported close to the diagonal $x=y$ is a convolution against an integrable function and
hence it maps $L^2$ to itself.
\end{proof}

If the map $VR_0^\pm(\lambda^2)$, or one of its derivatives (with respect to $\lambda$),
is applied enough times to a locally integrable function with fast decay, the result will be locally in $L^2$.
Any subsequent applications of the free resolvent will yield functions in weighted $L^2$ spaces.  Each time the Limiting
Absorption Principle is invoked, it improves the norm bounds by a factor of $\japanese[\lambda]^{-1+\eps}$
until eventually, some polynomial decay in $\lambda$ is achieved. Our primary estimate of this form is given below.
\begin{corollary} \label{enough}
Suppose $|V(x)| \leq C\japanese[x]^{-\beta}$ for some $\beta > n+3$.
Let $m_0 > \frac{n^2}2$ and $0 \leq j \leq \frac{n}2+1$.  Then
\begin{equation}\label{enough!!}
\Big\|\Big(\frac{d}{d\lambda}\Big)^j\big[VR_0^\pm(\lambda^2)\big]^{m_0}f
\Big\|_{L^{2,\sigma}} \lesssim
\lambda^{-j}\japanese[\lambda]^{j +1-2n}\norm[f][L^{1,\frac32}]
\end{equation}
for any $\sigma < \beta - (\frac{n+3}2)$.
\end{corollary}
\begin{proof}
The lower bound of $\frac{n^2}2$ is not intended to be sharp and was
obtained in the following manner:
It requires about $\frac{n}4$ iterations of $VR_0^\pm(\lambda^2)$ to smooth
an integrable function to local $L^2$ behavior (see Proposition~\ref{smoothing}) and one more to
reach a weighted $L^2$ space (see Proposition~\ref{weighted}).  Also, $\frac{n}{2}+1$ powers of
$VR_0^\pm(\lambda^2)$ can be lost to derivatives which we bound using Corollary~\ref{deriv in weighted}.
For each of these $\frac{3n+8}4$ operations, we have established only a crude bound which grows
like $\lambda^{\frac{n-3}2}$.
According to Lemma~\ref{lap}, each time the Limiting Absorption Principle is invoked,
this reduces the degree of polynomial growth by $1-\eps$,
so it needs to be done approximately $\frac{(3n+8)(n-3)}{8} + 2n-1 $ times.
Setting $m_0 > \frac{n^2}2$ is sufficient to obtain \eqref{enough!!}.
\end{proof}

We will also need the following mapping properties of $\Im R_0(\lambda^2)$.

\begin{proposition}\label{ImR0 prop}
Let $0\leq j\leq \frac{n}{2}+1$. Then, for $\sigma>\frac{n+3}{2}$ we have
\begin{align}\label{deriv ImR0 L2}
\Big\|\Big(\frac{d}{d\lambda}\Big)^j \Im R_0(\lambda^2)\Big\|_{L^{2,\sigma}\to L^{2,-\sigma}}
\lesssim \lambda^{n-2-j}\langle \lambda\rangle^{\frac{3}{2}}.
\end{align}
Moreover, assuming $|V(x)|\leq C\langle x\rangle^{-\beta}$ for some $\beta>n+3$, we have
\begin{align}\label{deriv ImR0 L1}
\Big\|\Big(\frac{d}{d\lambda}\Big)^j V\Im R_0(\lambda^2)\Big\|_{L^{1,\frac32}\to L^{1,\frac32}}
\lesssim \lambda^{n-2-j}\langle \lambda\rangle^{\frac{3}{2}}
\end{align}
while, for $m\geq 2m_0>n^2$, $\sigma>\frac{n+3}2$, and $\beta>2\sigma$, we have
\begin{align}\label{deriv Im VR0}
\Big\|\Big(\frac{d}{d\lambda}\Big)^j \Im [VR_0^+(\lambda^2)]^m f\Big\|_{L^{2,\sigma}}
\lesssim \lambda^{n-2-j}\langle \lambda\rangle^{j+\frac{5}{2}-2n +\frac{n^2(n-3)}{4}}\|f\|_{L^{1,\frac32}}.
\end{align}
\end{proposition}

\begin{proof}
From \eqref{ImR0}, we have the following formula for the kernel of $\Im R_0(\lambda^2)$:
\begin{align*}
\Im  R_{0}(\lambda^{2})(x,y)=a_{n-2,\frac{n-3}{2}}(\lambda|x-y|)\frac{e^{\pm i\lambda|x-y|}}{|x-y|^{n-2}}.
\end{align*}
Derivatives can affect $\Im R_0(\lambda^2)$ in two ways: Whenever a derivative falls on the symbol, this has the
effect of reducing the power of $\lambda$ by one. If a derivative falls on the phase, this has the effect of
increasing the power of $|x-y|$ by one. Hence, using the calculus of the symbols $a_{i,j}$, we get
\begin{align*}
\Big(\frac{d}{d\lambda}\Big)^j\Im R_0(\lambda^2)(x,y)
&=\sum_{\substack{j_1+j_2=j \\j_1,j_2\geq 0}}\lambda^{-j_1}a_{n-2,\frac{n-3}2}(\lambda|x-y|)
    \frac{e^{\pm i\lambda|x-y|}}{|x-y|^{n-2-j_2}}\\
&=\sum_{\substack{j_1+j_2=j \\j_1,j_2\geq 0}}\lambda^{n-2-j}a_{j_2,j_2-\frac{n-1}2}(\lambda|x-y|)
    e^{\pm i\lambda|x-y|}.
\end{align*}
Thus,
\begin{align}\label{deriv ImR0}
\Big|\Big(\frac{d}{d\lambda}\Big)^j\Im R_0(\lambda^2)(x,y)\Big|
\lesssim \lambda^{n-2-j}\langle \lambda|x-y|\rangle^{j-\frac{n-1}2}.
\end{align}

The estimate \eqref{deriv ImR0 L2} follows from \eqref{deriv ImR0} and
\begin{align}\label{showshow}
\int \frac{\langle x\rangle^{-2\sigma}\langle y\rangle^{-2\sigma}}{\langle \lambda|x-y|\rangle^{n-1-2j}}dxdy
\lesssim \langle \lambda\rangle^{3}.
\end{align}
For $0\leq j\leq \frac{n-1}{2}$, \eqref{showshow} follows from the bound
$\langle \lambda|x-y|\rangle^{-n+1+2j}\lesssim 1$; the resulting integral is finite whenever $\sigma>\frac{n}{2}$.
For $\frac{n-1}{2}<j\leq \frac{n}{2}+1$, we first apply Lemma~\ref{easylem} to the integral in the variable $y$ to obtain
\begin{align*}
\int \frac{\langle y\rangle^{-2\sigma}}{\langle \lambda|x-y|\rangle^{n-1-2j}}dy
\lesssim \lambda^{-n+1+2j}\int \frac{\langle y\rangle^{-2\sigma}}{|x-y|^{n-1-2j}}dxdy
\lesssim \langle \lambda\rangle^{3}\langle x\rangle^{-n+1+2j}.
\end{align*}
The remaining integral in the variable $x$ is finite under our assumptions on $\sigma$.

In view of \eqref{deriv ImR0}, the estimate \eqref{deriv ImR0 L1} follows from
\begin{align}\label{show}
\sup_{x\in\R^n} \langle x\rangle^{-\frac32}\int \frac{dy}{\langle y\rangle^{\beta-\frac32}\langle \lambda|x-y|\rangle^{\frac{n-1}2-j}}\lesssim1.
\end{align}
To see \eqref{show} one considers separately the cases $0\leq j\leq \frac{n-1}{2}$ and
$\frac{n-1}{2}<j\leq \frac{n}{2}+1$, bounding $\langle \lambda|x-y|\rangle^{-\frac{n-1}2+j}\lesssim 1$ in the former case and
applying Lemma~\ref{easylem} with $\sigma=\beta-\frac32$ and $\mu=\frac{n-1}{2}-j$ in the latter case.

We turn now to \eqref{deriv Im VR0}. We rewrite $\Im [VR_0^+(\lambda^2)]^m= [VR_0^+(\lambda^2)]^m- [VR_0^-(\lambda^2)]^m$
using the following algebraic identity:
\begin{equation} \label{algebraic identity}
\prod_{k=0}^M A_k^+  - \prod_{k=0}^M A_k^-  =
\sum_{l=0}^M \Big(\prod_{k=0}^{l-1} A_k^-\Big)\big(A_l^+ - A_l^-\big)
             \Big(\prod_{k=l+1}^M A_k^+\Big).
\end{equation}
Then,
\begin{align}\label{diff}
\Im [VR_0^+(\lambda^2)]^m
=\sum_{\substack{m_1+m_2=m \\m_1,m_2\geq 0}}[VR_0^-(\lambda^2)]^{m_1}V\Im R_0 [VR_0^+(\lambda^2)]^{m_2}.
\end{align}
We treat the cases $m_1<m_0$ and $m_2<m_0$ separately. In the first case, use Corollary~\ref{enough} for
$[VR_0^+(\lambda^2)]^{m_2}$, \eqref{deriv ImR0 L2}, and Corollary~\ref{deriv in weighted} for
$[VR_0^-(\lambda^2)]^{m_1}$ to derive the claim. In the second case, use the weighted $L^1$ bound in
Proposition~\ref{smoothing} for $[VR_0^+(\lambda^2)]^{m_2}$, \eqref{deriv ImR0 L1}, and Corollary~\ref{enough} for
$[VR_0^-(\lambda^2)]^{m_1}$ to obtain \eqref{deriv Im VR0}.
\end{proof}

We also record the following lemma whose proof is just an exercise
in integration by parts:
\begin{lemma}\label{intparts}
Given $a\in C^{\infty}_c(\R\setminus\{0\})$, we have
\begin{displaymath}
\bigl|\int_{\R}e^{i t\lambda^2}\lambda a(\lambda)d\lambda\bigr|
   \lesssim |t|^{-N} \sum_{s=0}^N \Bigl|\int_{\R}e^{i t\lambda^2}
\lambda^{s+1-2N} a^{(s)}(\lambda) d\lambda\Bigr|,
\end{displaymath}
for every $N\geq 0$.
\end{lemma}

%
%
%
%

\section{Dispersive Estimate for the Final Term} \label{sec:tail}

In this section we will show that the tail \eqref{tail} of the
finite Born series expansion \eqref{Born series} obeys dispersive
estimates for any potential $V$ satisfying $|V(x)|\lesssim \langle
x\rangle^{-\beta}$, provided we take $\beta$ and $m$ large
enough.

%
%

\begin{theorem}\label{dispersive tail}
Assume that the potential $V$ satisfies $|V(x)|\lesssim \langle
x\rangle^{-\beta}$ for some $\beta>\frac{3n+5}2$ and that
$m>n^2$. Then
\begin{align}
\sup_{x,y\in\R^n}\Bigl|\Im \int_{0}^{\infty}e^{it\lambda^{2}}\lambda
\bigl\{R_{0}^+(\lambda^{2})&V[R_{0}^+(\lambda^{2})V]^{m}
S^+(\lambda^2)R_{0}^+(\lambda^{2}) \notag \\
&\times \ [VR_{0}^+(\lambda^{2})]^{m}VR_{0}^+(\lambda^{2})\bigr\}(x,y)
d\lambda\Bigr|  \lesssim |t|^{-\frac{n}{2}}. \label{tail integral1}
\end{align}
\end{theorem}

\begin{remark} The condition $\beta > \frac{3n+5}2$ is not intended to be
 sharp.  Since the function we eventually construct as a counterexample
has compact support, decay conditions are not a matter of primary concern.
\end{remark}

There are numerous oscillatory components in this integral, which suggests
the use of stationary phase methods.  Although it appears natural to take the
critical point to be $\lambda = 0$, this turns out not to be the best choice.
Define the functions
$G_{\pm,x}(\lambda^{2})(\cdot) := e^{\mp i\lambda|x|}R_{0}^{\pm}(\lambda^{2})
(\cdot,x)$.  The expression in \eqref{dispersive tail} can be rewritten as
$I^+(t,x,y) - I^-(t,x,y)$, where
\begin{align} \label{Ipm}
I^\pm(t,x,y) &:= \int_{0}^{\infty}e^{it\lambda^{2}}e^{\pm i\lambda(|x|+|y|)}
\lambda \big\langle S^\pm(\lambda^2)R_{0}^\pm(\lambda^{2})
[VR_{0}^\pm(\lambda^{2})]^{m} VG_{\pm,y}(\lambda^{2}), \notag \\
& \hskip 2.5in [VR_{0}^\mp(\lambda^{2})]^{m}VG_{\mp,x}(\lambda^{2})
\big\rangle d\lambda \notag\\
& = \int_0^\infty e^{it\lambda^2}e^{\pm i\lambda(|x|+|y|)}
b^\pm_{x,y}(\lambda^2)\, d\lambda.
\end{align}
It suffices to show that $|I^+(t,x,y)-I^-(t,x,y)| \lesssim |t|^{-\frac{n}2}$
uniformly in $x$ and $y$.

The first step is to establish some properties (including existence)
of the operators $S^\pm(\lambda^2)$.  This is the crux of the Limiting
Absorption Principle for perturbed resolvents.  We sketch the details below.
\begin{proposition} \label{limabs}
Suppose $|V(x)| \leq C \japanese[x]^{-\beta}$ for some $\beta > \frac{n+1}{2}$
and also that zero energy is neither an eigenvalue nor a resonance of
$H = -\Delta + V$.  Then
$$\sup_{\lambda\geq0} \norm[S^\pm(\lambda^2)][L^{2,-\sigma}\to L^{2,-\sigma}]
 \ < \ \infty$$
for all $\sigma \in (\frac12, \beta-\frac12)$.
\end{proposition}
\begin{proof}
Under our assumptions, \eqref{deriv free res bound} and Proposition~\ref{impprop} imply that $R_0^\pm(\lambda^2)V$
is a compact operator on the space $L^{2,-\sigma}$.  The Fredholm alternative then guarantees the
existence of $S^\pm(\lambda^2)$ unless there exists a nonzero function $g \in L^{2,-\sigma}$ satisfying
$g = -R_0^\pm(\lambda^2)Vg$.

For $\lambda>0$, as $g = -R_0^\pm(\lambda^2)Vg$ is formally equivalent to $(-\Delta+V)g=\lambda^2g$,
it follows by a theorem of Agmon \cite{agmon} (see also \cite[Section XIII.8]{RSIII}) that
$g$ is in fact an eigenfunction, that is, $g\in L^2$. As positive imbedded eigenvalues do not exist by Kato's theorem
(see, for example, \cite[Section XIII.8]{RSIII}), we must have $g \equiv 0$.

When $\lambda = 0$, the free resolvent $R_0(0)$ is a scalar multiple of $I_2$. Since we are in dimension $n \geq 4$,
it is possible to improve the decay of $g$ by a bootstrap argument to obtain $g\in L^{2,-\sigma'}$ for all $\sigma'>0$;
in dimension $n\geq 5$, it is in fact possible to bootstrap all the way to $g\in L^2$.
In other words, zero energy would have to be either an eigenvalue or a resonance of $H$, contradicting our assumptions.
Thus, we must have $g\equiv 0$.

To obtain a uniform bound for $S^\pm(\lambda^2)$, note that by Lemma~\ref{lap} we have
$$
\|R_0^\pm(\lambda^2)V\|_{L^{2,-\sigma}} \lesssim \japanese[\lambda]^{-1+\eps}.
$$
Thus $I + R_0^\pm(\lambda^2)V$ converges to the identity as $\lambda \to\infty$.  Its inverse,
$S^\pm(\lambda^2)$, will thus have operator norm less than 2 for all $\lambda > \lambda_0$.
On the remaining interval, $\lambda \in [0,\lambda_0]$, observe that the family of operators
$R_0^\pm(\lambda^2)$ varies continuously with $\lambda$.  By continuity of inverses,
$S^\pm(\lambda^2)$ is continuous and bounded on this compact interval.
\end{proof}

Derivatives of $S^\pm(\lambda^2)$ can be taken using the identity
$$\frac{d}{d\lambda}S^\pm(\lambda^2) =
  - S^\pm(\lambda^2)\frac{d}{d\lambda}\big(R_0^\pm(\lambda^2)\big)
\,VS^\pm(\lambda^2).  $$
From this, Corollary \ref{deriv in weighted}, and Proposition~\ref{limabs}, it follows that for
$1\leq j\leq \frac{n}{2}+1$,
\begin{align}\label{S}
\Big\|\Big(\frac{d}{d\lambda}\Big)^jS^{\pm}(\lambda^2)\Big\|_{L^{2,-\sigma}\to L^{2,-\sigma}}
\lesssim \lambda^{-j}\langle \lambda\rangle^{j+\frac{n-3}{2}},
\end{align}
provided $\frac{1}{2}+j<\sigma<\beta-(\frac{1}{2}+j)$ and $\beta>\frac{n+1}{2}+j$. Moreover, it becomes clear that
$R_V^\pm(\lambda^2) = S^\pm(\lambda^2)R_0^\pm (\lambda^2)$ and its derivatives have mapping properties comparable to
those of the free resolvent.

We now have estimates for every object in \eqref{Ipm} except
for the functions $G_{\pm,y}(\lambda^2)$.  These follow from another
straightforward computation.
\begin{proposition} \label{Gboundprop}
Suppose $|V(x)| \leq C\japanese[x]^{-\beta}$ for some $\beta > \frac{3n+5}2$.
Then for each $0 \le j \le \frac{n}2+1$,
\begin{equation} \label{Gbound}
\Big\| V(\cdot)\Big(\frac{d}{d\lambda}\Big)^{j}G_{\pm,y}(\lambda^2)(\cdot)
  \Big\|_{L^{1,\frac32}} \lesssim  \frac{\lambda^{-j}}{\japanese[y]^{n-2}}
   + \frac{\lambda^{\frac{n-3}2-j}}{\japanese[y]^{\frac{n-1}2}}+
  \frac{\lambda^{\frac{n-3}2}}{\japanese[y]^{\frac{n-1}2}}.
\end{equation}
\end{proposition}
\begin{proof}
Write out the function $G_{\pm,y}(\lambda^2)$ in the form
$$
G_{\pm,y}(\lambda^2)(x) = a_{0,\frac{n-3}2}(\lambda|x-y|) \frac{e^{\pm i\lambda(|x-y|-y)}}{|x-y|^{n-2}}.
$$
Derivatives can affect $G_{\pm,y}$ in one of two ways.  Whenever a
derivative falls on the symbol, it has the effect of reducing the power of
$\lambda$ by one (this property was utilized previously in
Section~\ref{sec:lemmas}).  When derivatives fall on the exponential factor,
the effect is to multiply by $|x-y|-|y|$, which is smaller than $\japanese[x]$.
Thus, for $0 \le j \le \frac{n}2+1$,
\begin{equation*}
\Big(\frac{d}{d\lambda}\Big)^jG_{\pm,y}(\lambda^2)(x)
=\sum_{\substack{j_1+j_2=j \\j_1,j_2\geq 0}}\lambda^{-j_1}a_{0,\frac{n-3}2}(\lambda|x-y|)
    \frac{(|x-y|-|y|)^{j_2}}{|x-y|^{n-2}}e^{\pm i\lambda(|x-y|-y)}
\end{equation*}
and hence
\begin{equation*}
\Big|\Big(\frac{d}{d\lambda}\Big)^jG_{\pm,y}(\lambda^2)(x)\Big|
\lesssim \sum_{\substack{j_1+j_2=j \\j_1,j_2\geq 0}}
    \Big(\lambda^{-j_1}\frac{\langle x\rangle^{j_2}}{|x-y|^{n-2}}+\lambda^{\frac{n-3}{2}-j_1}\frac{\langle x\rangle^{j_2}}{|x-y|^{\frac{n-1}{2}}}\Big).
\end{equation*}
The result now follows from Lemma~\ref{easylem} provided $\beta>\frac{3n+5}{2}$.

\end{proof}

\begin{proof}[Proof of Theorem~\ref{dispersive tail}]
Consider first what happens if $|t| \leq 4$.  The bounds established in
Corollary~\ref{enough} (for $j=0$), Proposition~\ref{limabs}, and Proposition~\ref{Gboundprop} (for $j=0$)
show that the function $b^\pm_{x,y}(\lambda^2)$ in \eqref{Ipm}
is smaller than $\japanese[\lambda]^{-2}$ uniformly in $x$ and $y$.
This bounds the value of $I^+(t,x,y)-I^-(t,x,y)$ by a constant, which is
less than $|t|^{-\frac{n}2}$ as desired.

For the remainder of the calculation we will assume that $|t| > 4$.

Let $\rho: \R\to\R$ be a smooth even cutoff function which is identically one on
the interval $[-1, 1]$ and identically zero outside $[-2,2]$.
Let $b_{x,y,1}^\pm(\lambda^2)  := \rho(|t|^{\frac12}\lambda)b^\pm_{x,y}(\lambda^2)$
and $b_{x,y,2}^\pm := b_{x,y}^\pm - b_{x,y,1}^\pm$ and define $I^\pm_1(t,x,y)$, $I^\pm_2(t,x,y)$ accordingly.
For simplicity, the dependence on $x$ and $y$ will be
suppressed whenever possible.

We consider the integrals $I^\pm_2(t,x,y)$ first.

\noindent {\bf Case 1: $|x| + |y| \geq |t|$.}  At least one of $|x|$, $|y|$
is greater than $\frac{|t|}{2}$; without loss of generality assume it is $|y|$.
Then $|y|^{-1}\leq 2|t|^{-1}<|t|^{-\frac{1}{2}}$ and hence $|y|^{-1}$ does not belong to $\text{supp}\ b_2^{\pm}$.
Moreover, for $\lambda \geq |y|^{-1}$, Proposition~\ref{Gboundprop} yields the bound
\begin{align}\label{y}
 \norm[V\big({\textstyle\frac{d}{d\lambda}}\big)^jG_{\pm,y}(\lambda^2)][L^{1,\frac32}]
\lesssim \frac{\lambda^{\frac{n-3}2-j}\japanese[\lambda]^j}{\japanese[y]^{\frac{n-1}2}}
\lesssim |t|^{\frac{1-n}2} \lambda^{\frac{n-3}2-j}\japanese[\lambda]^j.
\end{align}
To bound $G_{\pm,x}(\lambda^2)$, we use
\begin{align}\label{x}
\norm[V(\tfrac{d}{d\lambda})^j  G_{\pm,x}(\lambda^2)][L^{1,\frac32}]
 \lesssim \lambda^{-j}   \japanese[\lambda]^{j + \frac{n-3}2}.
\end{align}
No additional improvement can be gained here, because the size of $|x|$ is unknown.

By Corollary \ref{deriv in weighted}, Corollary~\ref{enough}, Proposition~\ref{limabs}, \eqref{S}, \eqref{y}, and
\eqref{x}, we can deduce
$$
|b_2^\pm(\lambda^2)|
\lesssim |t|^{\frac{1-n}2}\lambda^{\frac{n-1}{2}}\langle \lambda\rangle^{-3n-1}
\lesssim |t|^{\frac{1-n}2}\langle \lambda\rangle^{-\frac{5n+3}{2}}
$$
and
$$
\Big|\frac{d}{d\lambda} b_2^\pm(\lambda^2)\Big|
\lesssim |t|^{\frac{1-n}2}\lambda^{\frac{n-3}2}\langle \lambda\rangle^{-\frac{5n+3}{2}}.
$$
Applying stationary phase around the critical
point $\lambda_0 =\mp \frac{|x|+|y|}{2t}$ and integrating by parts once away from the critical point, it follows that
$|I^\pm_2(t)| \lesssim |t|^{-\frac{n}2}$.

\noindent {\bf Case 2: $|t|^{\frac12} \leq |x|+|y| <|t|$}.
Again, assume without loss of generality that $|y| \geq \frac12|t|^{\frac12}$. Therefore, for
$\lambda \in \text{supp}\ b_2^{\pm}$ we have $|y|\geq \frac{1}{2}|\lambda|^{-1}$, which implies
$$
\norm[V\big({\textstyle\frac{d}{d\lambda}}\big)^jG_{\pm,y}(\lambda^2)][
 L^{1,\frac32}] \lesssim \frac{\lambda^{\frac{n-3}2-j}\japanese[\lambda]^j}{
 \japanese[y]^{\frac{n-1}2}}
$$
For $G_{\pm,x}(\lambda^2)$ we will use \eqref{x}.

The critical point for the phase occurs at $\lambda_0 = \mp\frac{|x|+|y|}{2t}$,
which is comparable in size to $\frac{|y|}{|t|}$ and greater than
$\frac12|t|^{-\frac12}$.  In the interval $[\lambda_0-\frac14|t|^{-\frac12},
\lambda_0+\frac14|t|^{-\frac12}]$ we have the size estimate
$$
|b^\pm_2(\lambda^2)| \ \sim \ \Big(\frac{|\lambda_0|}{|y|}\Big)^{\frac{n-1}2}  \ \sim \ |t|^{-\frac{n}2+\frac12}.
$$
An application of stationary phase yields the desired bound on this interval.

Away from the critical point, the derivatives of $b_2^\pm(\lambda^2)$ obey the following bounds
\begin{equation}\label{deriv b2}
\Big|\Big(\frac{d}{d\lambda}\Big)^jb_2^\pm(\lambda^2)\Big| \lesssim
\frac{\lambda^{\frac{n-1}2-j}\japanese[\lambda]^{j-\frac{5(n+1)}2}}{\japanese[y]^{\frac{n-1}2}}
\end{equation}
for all $0 \le j \le \frac{n}2+1$.

Over the intervals $[|t|^{-\frac12},\lambda_0 - \frac14|t|^{-\frac12}]$
and $[\lambda_0+\frac14|t|^{-\frac12}, 2\lambda_0]$, \eqref{deriv b2} becomes
\begin{align*}
\Big|\Big(\frac{d}{d\lambda}\Big)^jb_2^\pm(\lambda^2)\Big|
&\lesssim \frac{\lambda_0^{\frac{n-1}2}}{\japanese[y]^{\frac{n-1}2}}\lambda^{-j}\japanese[\lambda]^{j-\frac{5(n+1)}2}\\
&\lesssim |t|^{-\frac{n}{2}+\frac{1}{2}}\lambda^{-j}
\end{align*}
for all $0 \le j \le \frac{n}2+1$. As on this region $\lambda-\lambda_0\gtrsim |t|^{-\frac{1}{2}}$, each integration
by parts in \eqref{Ipm} gains us a factor of $|t|^{-\frac12}$. Thus, integrating by parts twice (i.e., taking $j=2$)
and recalling that in this case $\lambda_0\gtrsim |t|^{-\frac12}$, we obtain the desired dispersive estimate.

Over the interval $[2\lambda_0,1]$ (where $\lambda-\lambda_0\geq \frac12\lambda$), we use \eqref{deriv b2} and
the assumption $|y|\geq \frac12|t|^{\frac12}$ to get
\begin{align*}
\Big|\Big(\frac{d}{d\lambda}\Big)^jb_2^\pm(\lambda^2)\Big|
\lesssim \lambda^{\frac{n-1}2-j}|t|^{-\frac{n-1}4}
\end{align*}
for all $0 \le j \le \frac{n}2+1$. To obtain the desired decay in $t$, it is necessary to integrate by parts
at least $\frac{n+1}4$ times.

On the interval $[1,\infty]$, \eqref{deriv b2} implies that $b_2^\pm(\lambda^2)$ and its
derivatives all decay faster than $\japanese[\lambda]^{-2}  \japanese[y]^{\frac{1-n}2}$.
Using again the assumption $|y| \geq\frac12 |t|^{\frac12}$ and integrating by parts another $\frac{n+1}4$
times, we obtain the desired dispersive estimate.

\noindent {\bf Case 3: $|x|, |y| < |t|^{\frac12}$.}  This time,
the critical point $\lambda_0 = \mp\frac{|x|+|y|}{2t}$ lies outside the support
of $b^\pm_2(\lambda^2)$.  Therefore, one could safely integrate by parts;
however, the lack of a lower bound for $|x|$ and $|y|$ limits the usefulness
of estimates like \eqref{Gbound} in the regime $\lambda<1$. Without loss of generality, assume $|y| \geq |x|$.

For $\lambda \geq 1$, $b_2^\pm(\lambda^2)$ and its derivatives decay rapidly. Indeed, by
Corollary ~\ref{deriv in weighted}, Corollary~\ref{enough}, Proposition~\ref{limabs}, \eqref{S}, and
Proposition~\ref{Gboundprop}, for $\lambda\geq 1$ and $0\leq j\leq \frac{n}{2}+1$ we get
$$
\Big|\Big(\frac{d}{d\lambda}\Big)^j b_2^{\pm}(\lambda^2)\Big|
\lesssim \lambda^{-2n-3}\Big(\frac{1}{\langle x\rangle^{n-2}}+\frac{1}{\langle x\rangle^{\frac{n-1}{2}}}\Big)
                        \Big(\frac{1}{\langle y\rangle^{n-2}}+\frac{1}{\langle y\rangle^{\frac{n-1}{2}}}\Big).
$$
As the powers of $\japanese[x]$ and $\japanese[y]$ in the denominator may not make a meaningful contribution
(if $x,y$ are small), it is necessary to integrate by parts at least $\frac{n}2$ times in order to generate the
desired $|t|^{-\frac{n}2}$ decay or better.

The regime $\lambda \in [\japanese[y]^{-1},1]$ is similar to the interval $[2\lambda_0,1]$ in the previous case.
Indeed,
\begin{align*}
 \norm[V\big({\textstyle\frac{d}{d\lambda}}\big)^jG_{\pm,y}(\lambda^2)][L^{1,\frac32}]
\lesssim \frac{\lambda^{\frac{n-3}2-j}\japanese[\lambda]^j}{\japanese[y]^{\frac{n-1}2}}
\lesssim \frac{\lambda^{\frac{n-3}2-j}}{\japanese[y]^{\frac{n-1}2}}
\end{align*}
and
\begin{align*}
\norm[V(\tfrac{d}{d\lambda})^j  G_{\pm,x}(\lambda^2)][L^{1,\frac32}]
 \lesssim \lambda^{-j}   \japanese[\lambda]^{j + \frac{n-3}2}
 \lesssim \lambda^{-j}
\end{align*}
for all $0\leq j\leq \frac{n}{2}+1$. Thus,
\begin{equation*}
\Big|\Big(\frac{d}{d\lambda}\Big)^jb_2^\pm(\lambda^2)\Big|
\lesssim \frac{\lambda^{\frac{n-1}2-j}\japanese[\lambda]^{j-\frac{5(n+1)}2}}{\japanese[y]^{\frac{n-1}2}}
\lesssim \frac{\lambda^{\frac{n-1}2-j}}{\japanese[y]^{\frac{n-1}2}}
\end{equation*}
for all $0 \le j \le \frac{n}2+1$. Integrating by parts $\frac{n}{2}\leq N \leq\frac{n}2+1$ times is more than
enough to create polynomial decay in $\lambda$:
\begin{align*}
|t|^{-N} \int_{\langle y\rangle^{-1}}^1 (\lambda-\lambda_0)^{-N}\frac{\lambda^{\frac{n-1}2-N}}{\japanese[y]^{\frac{n-1}2}}d\lambda
\lesssim |t|^{-N}\japanese[y]^{2N-n}.
\end{align*}
Recalling that in this case we have $|y| <|t|^{\frac12}$, the resulting bound for this piece is
$|t|^{-\frac{n}{2}}$.

For the remaining interval, $ [|t|^{-\frac12},\japanese[y]^{-1}]$, we exploit instead the cancellation between
$R_0^+(\lambda^2)$ and $R_0^-(\lambda^2)$ using the algebraic identity \eqref{algebraic identity}.
We apply \eqref{algebraic identity} to $I_2^+(t,x,y)-I_2^-(t,x,y)$, where
\begin{align*}
I_2^\pm(t,x,y)
&=\int_0^\infty e^{it\lambda^2} \bigl(1-\rho(|t|^{\frac12}\lambda)\bigr)
  \lambda \bigl\{R_{0}^\pm(\lambda^{2})V[R_{0}^\pm(\lambda^{2})V]^{m}
S^\pm(\lambda^2)R_{0}^\pm(\lambda^{2})\\
& \hskip 2in \times \ [VR_{0}^\pm(\lambda^{2})]^{m}VR_{0}^\pm(\lambda^{2})\bigr\}(x,y)
\,d\lambda\\
&=\int_{0}^{\infty}e^{it\lambda^{2}}\bigl(1-\rho(|t|^{\frac12}\lambda)\bigr)
\lambda \big\langle \delta_y, R_{0}^\pm(\lambda^{2})V[R_{0}^\pm(\lambda^{2})V]^{m}
S^\pm(\lambda^2)R_{0}^\pm(\lambda^{2})  \\
& \hskip 2in \times \ [VR_{0}^\pm(\lambda^{2})]^{m}VR_{0}^\pm(\lambda^{2})\delta_x \big\rangle d\lambda\\
&=\int_{0}^{\infty}e^{it\lambda^{2}}c_{x,y,2}^{\pm}(\lambda^2)d\lambda.
\end{align*}
Each term in the resulting sum contains a factor of $R_0^+(\lambda^2)
-R_0^-(\lambda^2)$, an integral operator whose kernel is pointwise
dominated by $\lambda^{n-2}$ (see \eqref{deriv ImR0}).
This is even true if the cancellation falls on $S^+(\lambda^2)$ because we can write
$$
S^+(\lambda^2) - S^-(\lambda^2) = -S^-(\lambda^2)\big(R_0^+(\lambda^2)- R_0^-(\lambda^2)\big)VS^+(\lambda^2).
$$
We will integrate by parts $\frac{n+1}2$ times if $n$ is odd and $\frac{n}2+1$
times if $n$ is even. Our analysis relies on the estimates of Proposition~\ref{ImR0 prop}.

In place of the weighted $L^1$ estimate \eqref{Gbound}, we use the following two bounds for the two possible
initial functions on which the resolvents act. For $0\leq j\leq \frac{n}2+1$, we have
\begin{align}
\Big\|V(\cdot) \Big(\frac{d}{d\lambda}\Big)^j R_0^\pm(\lambda^2){\textstyle (\cdot,y)}\Big\|_{L^{1,\frac32}}
&\lesssim\frac{\lambda^{-j}}{\japanese[y]^{n-2}}, \label{show1}\\
\Big\|V(\cdot) \Big(\frac{d}{d\lambda}\Big)^{j}\big(\Im R_0(\lambda^2)\big)(\cdot,y)\Big\|_{L^{1,\frac32}}
&\lesssim\lambda^{n-2-j}. \label{show2}
\end{align}
To see \eqref{show1}, we use the pointwise bound
$$
\Big|\Big(\frac{d}{d\lambda}\Big)^jR_0^\pm(\lambda^2)(x,y)\Big|
\lesssim \lambda^{-j}I_2 +\lambda^{\frac{n-3}2}I_{\frac{n+1}2+j}
$$
and apply Lemma~\ref{easylem} to obtain
\begin{align*}
\Big\|V(\cdot) \Big(\frac{d}{d\lambda}\Big)^j R_0^\pm(\lambda^2){\textstyle (\cdot,y)}\Big\|_{L^{1,\frac32}}
\lesssim \frac{\lambda^{-j}}{\japanese[y]^{n-2}}+\frac{\lambda^{\frac{n-3}2}}{\japanese[y]^{\frac{n-1}2-j}}
\lesssim \frac{\lambda^{-j}}{\japanese[y]^{n-2}},
\end{align*}
where the last inequality holds for $\lambda \leq \japanese[y]^{-1}$.

Similarly, to prove \eqref{show2} we use \eqref{deriv ImR0}; applying Lemma~\ref{easylem} and
treating the cases $0\leq j\leq \frac{n-1}2$ and $\frac{n-1}2<j\leq \frac{n}{2}+1$ separately, we obtain
\begin{align*}
\Big\|V(\cdot) \Big(\frac{d}{d\lambda}\Big)^{j}\big(\Im R_0(\lambda^2)\big)(\cdot,y)\Big\|_{L^{1,\frac32}}
\lesssim \lambda^{n-2-j} \big(1 +\lambda^{\frac32}\langle y\rangle^{\frac{3}{2}}\big)
\lesssim\lambda^{n-2-j},
\end{align*}
again, for $\lambda \leq \japanese[y]^{-1}$.

Using the estimates in Proposition~\ref{ImR0 prop}, \eqref{show1}, and \eqref{show2}, we get
$$
\Big|\Big(\frac{d}{d\lambda}\Big)^{j}\big(c_{x,y,2}^+(\lambda^2)-c_{x,y,2}^-(\lambda^2)\big)\Big|\lesssim \frac{\lambda^{n-1-j}}{\langle y\rangle^{n-2}}.
$$
Thus, an application of Lemma~\ref{intparts} with $N=\frac{n+1}2$ for $n$ odd, or $N=\frac{n}{2}+1$ for $n$ even
yields the bound
\begin{align*}
|t|^{N}\int_{|t|^{-\frac12}}^{\langle y\rangle^{-1}} \frac {\lambda^{n-1-2N}}{\langle y\rangle^{n-2}}d\lambda
\lesssim |t|^{-\frac{n}2}.
\end{align*}

In each of the three cases discussed above, the difference $I^+_2(t,x,y)-I^-_2(t,x,y)$ is
seen to be smaller than $|t|^{-\frac{n}2}$.

To complete the proof of the theorem, we need to show
$$|I^+_1(t,x,y)-I^-_1(t,x,y)|\lesssim |t|^{-\frac{n}2} \ \ \text{for} \ \ |t|>4.$$
Here,
\begin{align*}
I_1^\pm(t,x,y)
&=\int_0^\infty e^{it\lambda^2} \rho(|t|^{\frac12}\lambda)
  \lambda \bigl\{R_{0}^\pm(\lambda^{2})V[R_{0}^\pm(\lambda^{2})V]^{m}
S^\pm(\lambda^2)R_{0}^\pm(\lambda^{2})\\
& \hskip 2in \times \ [VR_{0}^\pm(\lambda^{2})]^{m}VR_{0}^\pm(\lambda^{2})\bigr\}(x,y)
\,d\lambda\\
&=\int_{0}^{\infty}e^{it\lambda^{2}}\rho(|t|^{\frac12}\lambda)
\lambda \big\langle \delta_y, R_{0}^\pm(\lambda^{2})V[R_{0}^\pm(\lambda^{2})V]^{m}
S^\pm(\lambda^2)R_{0}^\pm(\lambda^{2})  \\
& \hskip 2in \times \ [VR_{0}^\pm(\lambda^{2})]^{m}VR_{0}^\pm(\lambda^{2})\delta_x \big\rangle d\lambda\\
&=\int_{0}^{\infty}e^{it\lambda^{2}}c_{x,y,1}^{\pm}(\lambda^2)d\lambda.
\end{align*}
Arguing as in Case 3 above, we see that
$$
\big|c_{x,y,1}^+(\lambda^2)-c_{x,y,1}^-(\lambda^2)\big|\lesssim \lambda^{n-1}.
$$
Thus,
$$
|I^+_1(t,x,y)-I^-_1(t,x,y)|\ \lesssim\ \int_0^{|t|^{-\frac12}} \lambda^{n-1}\,d\lambda \ \lesssim \ |t|^{-\frac{n}2}.
$$
This concludes the proof of Theorem~\ref{dispersive tail}.
\end{proof}

\section{Nondispersive Estimates}
\subsection{Nondispersive estimate for the term $l=1$} \label{sec:l=1}
To summarize the progress up to this point, we have decomposed the
perturbed resolvent $R_V^\pm(\lambda^2)$ into a finite Born series
with initial terms given by \eqref{first terms} and a tail given
by \eqref{tail}. In the previous sections, the contribution of the
tail was shown to satisfy a dispersive estimate at both high and
low energies.  The dispersive behavior of the full evolution
$e^{itH}P_{ac}(H)$ is therefore dictated by the contribution from
the initial terms of the Born series.

We show that there are potentials in the class
\[X=\bigg\{V\in C^\alpha(\R^n),\ \alpha<\tfrac{n-3}2, \ \text{supp}V
\subset B(0,5)\setminus B(0,\tfrac52)\bigg\}   \] that do not yield a dispersive estimate for the term
corresponding to $l=1$ in \eqref{first terms}. It will follow, via an argument in the next subsection, that
the entire expression \eqref{first terms} cannot satisfy a dispersive estimate either. To define the
class of potentials more precisely, let $X$ be the completion of the appropriately supported $C^\infty$
functions with respect to the $W^{\alpha,\infty}$-norm,
\[ \norm[f][X] := \norm[(1+\Delta)^{\alpha/2}f][\infty].   \]

Fix the points $x_0, y_0 \in \R^n$ so that $x_0$ is the unit vector in the
first coordinate direction and $y_0 = -x_0$.
Now let $f^{\eps}$ and $g^{\eps}$ be smooth approximations of $f=\delta_{x_0}$
and $g=\delta_{y_0}$ which are supported in $B(x_0,\eps)$ and $B(y_0,\eps)$,
respectively, and have unit $L^1$-norm.
Define the expression
\begin{align*}
a_1^L(t,\eps,V) :=
&t^{\frac{n}2}\int e^{it\lambda}\psi_L(\lambda) \langle
  [R_{0}^+(\lambda)(x,x_{1})V(x_1)R_{0}^+(\lambda)(x_{1},y) \\
&\hskip 0.5in - R_0^-(\lambda)(x,x_1)V(x_1)R_0^-(\lambda)(x_1,y)]
  f^\eps(x),g^\eps(y)\rangle\, dx dx_1 dy d \lambda \\
= &t^{\frac{n}2} \int I_{L}(t, |x-x_{1}|,|y-x_{1}|)V(x_{1})
   f^\eps(x)g^\eps(y)dx_{1}dxdy
\end{align*}
where $\psi$ can be any Schwartz function with $\psi(0) = 1$ and $\psi_L(\lambda) = \psi(\lambda/L)$.
Fubini's theorem is used to perform the $d\lambda$ integral first, noting that since $f^\eps$, $g^\eps$,
and $V$ all have disjoint support, the singularities of $R^{\pm}_0(\lambda)(x,x_1)$ and
$R^{\pm}_0(\lambda)(x_1,y)$ can be disregarded.

If the term corresponding to $l=1$ in the Born series \eqref{first terms}
satisfied a dispersive estimate, it would yield the bound
\begin{equation} \label{eq:a1dispersive}
\lim_{L\to \infty} |a_1^L(t,\eps,V)|
 \ \leq \ C(V) \norm[f^\eps][1]\norm[g^\eps][1] \ =\ C(V).
\end{equation}
Observe that $a_1^L(t,\eps,V)$ is linear in the last entry and can therefore be viewed as a family of
linear maps indexed by the remaining parameters $(L,t,\eps)$.  By the Uniform Boundedness Principle, if
a dispersive estimate for the $l=1$ term held for every potential $V \in X$, it would imply the sharper
inequality
\begin{equation}\label{a1bdd}
\sup_{L\geq 1} |a_1^L(t,\eps,V)| \leq C\norm[V][X].
\end{equation}
For $t \ll 1$ this will not be possible, thanks to the asymptotic description of the function $I_L(t,
|x-x_1|,|y-x_1|)$ stated below.

\begin{lemma} \label{lem:asymptotic}
Suppose $n \geq 3$ and $0 < t \leq 1$.  Let $\psi:\R\to\R$ be a Schwartz function with Fourier transform
supported in the unit interval and satisfying $\psi(0) = 1$, and $K$ a compact subset of $(0,\infty)$.
There exist constants $C_1, C_2 < \infty$ depending on $n$, $\psi$, and $K$ such that
\begin{align} \label{eq:asymptotic}
\Big| I_L(t,|x-x_1|,|x_1-y|)
 &- {\textstyle\frac{i}{2(-4\pi i\,t)^{n-\frac32}}
\Big(\frac{(|x-x_1|+|x_1-y|)^{n-2}}{|x-x_1|^{\frac{n-1}2} |x_1-y|^{\frac{n-1}{2}}}\Big)
e^{-i\frac{(|x-x_1|+|x_1-y|)^2}{4t}}} \Big| \notag\\
&\leq C_1 t^{-(n-\frac52)}
\end{align}
for all $L > C_2t^{-3}$ and $|x-x_1|, |x_1-y| \in K$.
If $t$ is held fixed, then the
remainder converges as $L\to\infty$ to a function $G(|x-x_1|,|x_1-y|,t)$
uniformly over all pairs of distances $|x-x_1|,|x_1-y| \in K$.
\end{lemma}
The proof of Lemma \ref{lem:asymptotic} is technical and is given in Section~6 below. An immediate
consequence of this lemma is the following

\begin{corollary}
Let $n \geq 3$, $0 <t \leq 1$, and $\eps < \frac12$.  The following bound is valid for all functions $V
\in X$ with $\norm[V][X] \leq 1$:
\begin{equation} \label{eq:a1infty}
\begin{aligned}
 \lim_{L\to\infty} &\bigg|a_1^L(t,\eps,V)
-  \frac{i\,t^{\frac{3-n}2}}{2(-4\pi i)^{n-\frac32}} \int
\bigg(\frac{(|x-x_1|+|x_1-y|)^{n-2}}{|x-x_1|^{\frac{n-1}2}
|x_1-y|^{\frac{n-1}{2}}} \bigg) \\
& \hskip 1.8in \times e^{-i\frac{(|x-x_1|+|x_1-y|)^2}{4t}}
    V(x_1) f^\eps(x) g^\eps(y)\, dx_1 dx dy \bigg| \\
&\leq Ct^{\frac{5-n}2}\norm[f^\eps][1]\norm[g^\eps][1].
\end{aligned}
\end{equation}
\end{corollary}

\begin{proof}
If $\eps <\frac12$, then we have $|x-x_1|,|x_1-y| \in [1,10]$
for every combination of points with $x\in \text{supp}(f^\eps)$, $y\in \text{supp}(g^\eps)$, $x_1 \in
\text{supp}(V)$. Thus the conditions of Lemma~\ref{lem:asymptotic} are satisfied, with the conclusion
that $I_L(t,\cdot,\cdot)$ converges uniformly as $L\to \infty$ to a bounded function in $x,x_1,y$.

The result then follows from the dominated convergence theorem and the observation that $\norm[V][1]
\leq C\norm[V][X] \leq C$ .
\end{proof}

If the integral in \eqref{eq:a1infty} were taken in absolute values, the resulting bound on
$a_1^L(t,\eps,V)$ would be of size $|t|^{\frac{3-n}2}$.  In dimension $n\geq 4$, this contrasts with the
desired estimate
$$
\lim_{L\to\infty}|a_{1}^{L}(t,\eps,V)|\leq C,
$$
which is uniform in $t$.  Furthermore, for a fixed small time $t$ it is not difficult to construct a
potential $V_t \in X$ which negates the oscillatory factor of $e^{-i(|x-x_1|+|x_1-y|)^2/(4t)}$.

Let $\phi$ be a smooth cutoff which is supported in the interval $[6,8]$
and $F:\R\to\R$ a nonnegative smooth function which satisfies
$F(s) = 0$ for all $s \leq 0$ and $F(s)=s$ for all $s \geq \frac12$.
Given a time $0<t\leq 1$, define
\begin{equation}
V_t(x_1) = C_n t^{\alpha}\phi(|x_0-x_1|+|x_1-y_0|)
    F\Big(\cos\Big(\frac{(|x_0-x_1|+|x_1-y_0|)^2}{4t}\Big)\Big).
\end{equation}
The constant $C_n$ will be chosen momentarily. It is perhaps unnecessary to modify the cosine function
with $F$; however, the positivity of $F$ does guarantee that zero energy will be neither an eigenvalue
nor a resonance of $-\Delta + V_t$.

\begin{proposition}
There exists a constant $C_n >0$ so that the function
$V_t$ defined above satisfies $\norm[V_t][X] \leq 1$ for all $0<t\leq  1$.
\end{proposition}
\begin{proof}
It is equivalent to show that in the absence of the coefficient $C_n$,
$\norm[V_t][X]$ would be bounded by a finite constant uniformly in $t$.

The support of $\phi(|x_0-x_1|+|x_1-y_0|)$ is located within an annular
region bounded by the ellipsoids with foci $x_0,y_0$ and major axes of
length $6$ and $8$, respectively.  As this region is bounded away from both
$x_0$ and $y_0$, the length sum
$|x_0-x_1| + |x_1-y_0|$ is a scalar $C^\infty$-function of $x_1$.

It follows that any sufficiently smooth function of
$\frac{|x_0-x_1|+|x_1-y_0|}{4t}$ on this domain should have
$C^\alpha$-norm controlled by $(1+t^{-\alpha})$. The leading
coefficient $t^{\alpha}$ then ensures that the $X$-norm will be
controlled by a uniform constant for all $|t| \leq 1$.  Finally,
multiplication by the fixed smooth cutoff
$\phi(|x_0-x_1|+|x_1-y_0|)$ only increases the norm by another
finite constant.
\end{proof}

Now it is a simple matter to show that $V_t$ produces a counterexample to \eqref{a1bdd} and hence to
\eqref{eq:a1dispersive} for $0<t\ll 1$.
\begin{proposition} \label{prop:a1bound}
Suppose $n > 3$.  There exist constants $T, C_1, C_2 > 0$ such that
if $0 < t \leq T$ and $0 < \eps < C_1t$, then
\begin{equation*}
\lim_{L\to\infty} \Big|a_1^L(t,\eps,V_t)\Big| \ge C_2
t^{-(\frac{n-3-2\alpha}2)}.
\end{equation*}
\end{proposition}
\begin{proof}
Start with the asymptotic integral formula in \eqref{eq:a1infty}.
For any choice of points $x\in \text{supp}(f^\eps)$,
$y\in \text{supp}(g^\eps)$, $x_1\in\text{supp}(V_t)$,
the expression $\frac{(|x-x_1|+|x_1-y|)^{n-2}}{(|x-x_1|\,|x_1-y|)^{(n-1)/2}}$
is a smooth positive function of size comparable to 1.

Consider what happens to the integral over $dx_1$ in the special case when $x = x_0$, $y = y_0$.  Then,
the oscillatory part of $V_t(x_1)$ is synchronized with the real part of
$e^{-i(|x-x_1|+|x_1-y|)^2/(4t)}$ so that the real part of the product is always positive and of size
approximately 1 on a set of approximately unit measure.  The real part of the integral is then bounded
below by a positive constant.

For arbitrary $x\in \text{supp}(f^\eps)$ and $y\in \text{supp}(g^\eps)$, it is possible to differentiate
under the integral sign in either of the variables $x$ or $y$ and each partial derivative is controlled
by $t^{-1}$. Thus the lower bound on the real part of the integral remains valid so long as
$|x-x_0|,|y-y_0| \lesssim t$, which is ensured by setting $\eps < C_1t$.

The definition of $V_t$ also includes a factor of $t^\alpha$. When this is substituted into
\eqref{eq:a1infty}, the resulting leading coefficient is proportional to $t^{-(\frac{n-3-2\alpha}2)}$.
There is also an error term of unknown sign, but with size controlled by $t^{-(\frac{n-5-2\alpha}2)}$.
This can be absorbed into the lower bound for any $0<t\leq T$, provided $T$ is chosen sufficiently
small.
\end{proof}

%
%
%
%

\subsection{Nondispersive Estimate for the Full Evolution}

\begin{theorem}\label{T}
Suppose $n > 3$.  There cannot exist a bound of the form
\begin{equation*}
\norm[e^{itH}P_{ac}f][\infty] \le C(V)|t|^{-\frac{n}2}\norm[f][1]
\end{equation*}
with $C(V) <\infty$ for every potential $V\in X$, $\norm[V][X] \leq 1$.
\end{theorem}
\begin{proof}
Assume the contrary and write $V = \theta W$ with $\norm[W][X] \leq 1$ and
$\theta \in [0,1]$.  By assumption, we would then have the bound
\begin{align}
|\langle e^{itH}P_{ac}f,g\rangle| &=
\frac{1}{2\pi} \sup_{L\geq 1} \Big|\int_0^\infty e^{it\lambda}
 \psi_L(\lambda)\langle [R_{\theta W}(\lambda+i0)-R_{\theta W}(\lambda-i0)]f,
  g\rangle\,d\lambda \Big| \notag \\
&\leq C(W,\theta)|t|^{-\frac{n}2} \norm[f][1]\norm[g][1] \label{Stone}
\end{align}
for $\psi$ as in Lemma~\ref{lem:asymptotic} and for every $f,g\in L^1\cap L^2$ and, in particular, for
the functions $f^\eps, g^\eps$ defined in subsection~\ref{sec:l=1}.

The finite Born series expansion \eqref{Born series} allows us to write the perturbed resolvent
$R_{\theta W}(\lambda\pm i0)$ as the sum of a polynomial of degree $2m+1$ in $\theta$ and a tail.  When
this is substituted into \eqref{Stone} above, along with the functions $f^\eps, g^\eps \in L^1\cap L^2$,
the tail is shown in Theorem~\ref{dispersive tail}
to be controlled by
$C|t|^{-\frac{n}{2}}\norm[f^\eps][1]\norm[g^\eps][1]$ for some $C$.
It follows that the initial terms must obey a similar bound. Write
this as
$$
\sup_{L\geq 1} \Big|P^L(\theta)\Big| := \sup_{L\geq 1} \Bigl|\sum_{k=0}^{2m+1}\theta^{k}a_{k}^{L}\Bigr|
\leq C(W,\theta),
$$
where the coefficients $\{a_k^L\}_{k=0}^{2m+1}$ of the polynomial $P^L$ are defined for each
$k\in\{0,1,..., 2m+1\}$ and $L < \infty$ by the formula
$$
a_{k}^{L}(t,\eps,W)=t^{\frac{n}2}\int e^{it\lambda}\psi_L(\lambda) \langle [R_{0}^+(\lambda)([-W
R_{0}^+(\lambda)]^{k} -R_0^-(\lambda)[-W R_0^-(\lambda)]^k]f^{\eps},g^{\eps}\rangle d\lambda.
$$

Denote by $\mathbb{V}$ the $2m+2$-dimensional space of all polynomials of degree $2m+1$, and consider
the linear maps from $\mathbb{V}$ into $\R^{2m+2}$ defined by
$$
P=\sum_{k=0}^{2m+1} a_{k}\theta^{k} \mapsto \{a_{0}, ..., a_{2m+1}\}
$$
and
$$
P=\sum_{k=0}^{2m+1} a_{k}\theta^{k} \mapsto \Bigl\{P\big(0\big),
 P\bigl(\tfrac{1}{2m+1}\bigr), \ldots, P\bigl(\tfrac{2m+1}{2m+1}\bigr)\Bigr\}.
$$

Clearly the two maps are bijections and thus one can express each coefficient $a_{k}$ as a linear
combination of the values $P(0),P(\frac{1}{2m+1}),..., P(\frac{2m+1}{2m+1})$. From our assumption that
$C(W,\theta) < \infty$ for every $0\leq \theta \leq 1$, it follows that each of the expressions
$|P^L(0)|, |P^L(\tfrac{1}{2m+1})|,\ldots, |P^L(\frac{2m+1}{2m+1})|$, as well as their maximum, is
bounded uniformly in $L\geq 1$.  One concludes that
\[  \sup_{L\geq 1} |a_1^L(t,\eps,W)| \leq C(W) <\infty  \]
for every $W \in X$ with $\norm[W][X] \leq 1$.

This, however, is precisely the same statement as \eqref{eq:a1dispersive}
which was already shown to be false.
\end{proof}

\section{Proof of Lemma \ref{lem:asymptotic}}
The main ingredients of Lemma~\ref{lem:asymptotic} are a recurrence relation
(in $n$) for the resolvent kernels
and explicit computations in dimensions 2 and 3.  With some abuse of notation,
define $R_n^\pm(\lambda)$  to be the free resolvent $\lim_{\eps\downarrow 0}
(-\Delta -\lambda \pm i\eps)^{-1}$ in $\R^n$.  The Stone formula dictates that
\begin{equation*} \begin{aligned}
\frac{1}{2\pi i} \int_0^{\infty} e^{it\lambda}
\big\langle[R_n^+(\lambda)-R_n^-(\lambda)]f, g\big\rangle \, d\lambda
 &= \big\langle e^{-it\Delta}f, g \big\rangle \\
 &= (-4\pi i\,t)^{-\frac{n}2} \iint_{\R^{2n}} e^{\frac{-i |x-y|^2}{4t}} f(x)
\bar{g}(y)\, dxdy
\end{aligned} \end{equation*}
for all $t \not= 0$ and $f, g$ (say) Schwartz functions.

Recall that the resolvents $R_n(z) = (-\Delta - z)^{-1}$ can be
defined for all $z \in \C \setminus \R^+$, and that
$R_n^\pm(\lambda)$ are the analytic continuations onto the
boundary from above and below, respectively. It follows that both
$R_n^+(\lambda)$ and $R_n^-(\lambda)$ can be defined for negative
values of $\lambda$.  Moreover, $[R_n^+(\lambda) - R_n^-(\lambda)]
= 0$ for all $\lambda \le 0$.  The integral above may therefore be
taken over the entire real line.

One further observation is that since $R_n^+(\lambda)$ is a holomorphic family of operators for
$\lambda$ in the upper halfplane and is uniformly bounded (as operators on $L^2$, for example) away from
the real axis, its inverse Fourier transform must be supported on the halfline $\{t \leq 0\}$.
Similarly, $R_n^-(\lambda)$, which is holomophic in the lower halfplane, has inverse Fourier transform
supported in $\{t \geq 0\}$. This leads to the conclusion
\begin{equation} \label{eq:Rtransform}
\int_\R e^{it\lambda} R_n^-(\lambda,|x|)\, d\lambda =
\left\{ \begin{aligned}
-2\pi i(-4\pi i\,&t)^{-\frac{n}{2}} e^{\frac{-i|x|^2}{4t}},&  &{\rm if}\ t>0 \\
    &0, & &{\rm if}\ t < 0\\
\end{aligned} \right.
\end{equation}
for all $x \in \R^n$.
Setting $|x| = r$ in the preceding identity leads to the recurrence relation
\begin{equation}
R_{n+2}^-(\lambda,r) = -\frac1{2\pi r} \ptl[][r] \big[R_n^-(\lambda,r)\big].
\end{equation}
The same identity also holds for $R_n^+(\lambda,r)$.

\subsection{The cases $n = 2,3$}

It should first be noted that the integral $\int_\R e^{it\lambda}R_n^-(\lambda,r)\, d\lambda$ in
\eqref{eq:Rtransform} is never absolutely convergent and is properly interpreted as the Fourier
transform of a distribution. As such, its behavior at $t=0$ requires additional clarification.

\begin{lemma}
For any fixed $r > 0$ and $n = 2,3$, the expression $\int_\R e^{it\lambda} R_n^-(\lambda,r)\, d\lambda$
agrees with the distribution $f$ given by
\begin{equation} \label{eq:Rdistribution}
(f,\phi) = -2\pi i(-4\pi i)^{-\frac{n}{2}}
\lim_{a\downarrow 0} \int_a^\infty
   t^{-\frac{n}2} e^{\frac{-i r^2}{4t}} \phi(t)\, dt
\end{equation}
for all Schwartz functions $\phi$.
\end{lemma}

\begin{proof}
Because of analyticity considerations, the identity above must be correct modulo distributions supported
on $t=0$. Let $\phi \in C^\infty_c(\R)$ have nonvanishing derivatives at $t=0$ and consider pairings of
the form $\langle f,N\phi(N\,\cdot)\rangle$.

On one hand, the function $t^{-n/2} e^{-i r^2/(4t)}\chi_{(0,\infty)}$ has a
continuous anti-derivative $I(t)$ with asymptotic behavior
$I(t) =  O(t^{2-\frac{n}2})$ as $t$ approaches zero.  Integrating by parts,
\begin{equation*}
\begin{aligned}
\lim_{a\downarrow 0} \int_a^\infty N \phi(Nt) t^{-\frac{n}2}
  e^{\frac{-i r^2}{4t}}\, dt &= - N\lim_{a\downarrow 0} \phi(Na)I(a)
     - \lim_{a\downarrow 0} \int_a^\infty N^2\phi'(Nt) I(t)\, dt \\
 &= -N^2\int_0^\infty \phi'(Nt) I(t)\, dt \\
 &= O(N^{\frac{n-2}{2}}) \quad {\rm in\ the\ limit\ }\ N \to \infty.
\end{aligned}
\end{equation*}

Meanwhile, the pairing $\langle f,N\phi(N\,\cdot)\rangle$ is defined by
Parseval's identity to be
\begin{equation*}
\big\langle f, N\phi(N(\cdot))\big\rangle =  \int_\R
   R_n^-(\lambda,r) \hat{\phi}(\lambda/N)\, d\lambda.
\end{equation*}
For fixed $r > 0$, the resolvent $R_n^-(\lambda,r)$
possesses the asymptotic expansion
\begin{equation*}
R_n^-(\lambda,r) = c_1r^{\frac{1-n}2}
\lambda^{\frac{n-3}{4}}e^{-ir\sqrt{\lambda}}
                                  + O(\lambda^{\frac{n-5}{4}})
\end{equation*}
as $\lambda\to\infty$ and is integrable near $\lambda=0$.  Thus, it has a continuous anti-derivative
$J(\lambda,r)$ which grows no faster than $O(\lambda^{\frac{n-1}{4}})$.  Integrating by parts,
\begin{equation*}
\begin{aligned}
\big\langle f, N\phi(N\,\cdot)\big\rangle &= - N^{-1} \int_\R
    J(\lambda,r) \hat{\phi}'(\lambda/N)\, d\lambda \\
&= O(N^{\frac{n-1}{4}}).
\end{aligned}
\end{equation*}

As $n=2,3$, the difference between the left and right sides of \eqref{eq:Rdistribution} grows no faster
than $O(N^{\frac12})$ when applied to the test functions $N\phi(Nt)$.  It is well-known that any nonzero
distribution $g$ supported on $t=0$ has the form $(g,\phi) = \sum_{k=1}^M c_k \phi^{(k)}(0)$, and would
therefore grow at least as fast as $O(N)$ when applied to the same family of test functions.
\end{proof}

Having established the inverse Fourier transform of $R_n^-(\lambda,r)$ for
each $r > 0$, it is possible to calculate the inverse Fourier transform of
any product $R_n^-(\lambda,r)R_n^-(\lambda,s)$ by taking convolutions.
Given a choice of $r,s,t > 0$,
\begin{equation} \label{eq:convolution}
\int_\R e^{it\lambda}R_n^-(\lambda,r)R_n^-(\lambda,s)\,d\lambda
 = \frac{-2\pi}{(-4\pi i)^n}
\int_0^t e^{-i(\frac{r^2}{4u} + \frac{s^2}{4(t-u)})}
  \frac{du}{u^{\frac{n}2}(t-u)^{\frac{n}2}}
\end{equation}
where the Fourier transform has introduced a normalizing factor of
$(2\pi)^{-1}$.  To make the complex exponential more manageable,
change variables to
\begin{equation*}
\frac{v}{t} = \frac{r^2}{u} + \frac{s^2}{t-u} - \frac{r^2+s^2}{t}
      = \frac{(t-u)^2r^2 + u^2s^2}{u(t-u)t}.
\end{equation*}
The range of possible values for $v$ is $[2rs, \infty)$.  Based on the
quadratic relationship
\begin{equation*}
(r^2+s^2+v)u^2 -(2r^2+v)tu + r^2t^2 = 0,
\end{equation*}
the variable substitutions for $u$ and $(t-u)$ are given by
\begin{equation*}
u = \bigg(\frac{2r^2}{2r^2+v \mp\sqrt{v^2-4r^2s^2}}\bigg)t, \qquad
t-u 
 = \bigg(\frac{\frac12\big(\sqrt{v+2rs}\mp\sqrt{v-2rs}\big)^2}{2r^2+v\mp
  \sqrt{v^2-4r^2s^2}}\bigg)t.
\end{equation*}
The substitution formula for the differentials is
\begin{equation*}
du = \pm t\bigg(\frac{r \big(\sqrt{v+2rs}\mp\sqrt{v-2rs}\big)}
 {2r^2+v \mp\sqrt{v^2-4r^2s^2}}\bigg)^2\frac{dv}{\sqrt{v^2-4r^2s^2}}.
\end{equation*}
Making all appropriate substitutions and correctly accounting for the fact that each value of $v > 2rs$
is attained twice in $u\in(0,t)$, the integral in \eqref{eq:convolution} becomes
\begin{align}
\int_\R e^{it\lambda}R_2^-(\lambda,r)R_2^-(\lambda,s)\,d\lambda
&=  \frac{1}{4\pi t} e^{-i\frac{r^2+s^2}{4t}}
   \int_{2rs}^\infty \frac{e^{-i\frac{v}{4t}}}{\sqrt{v^2-4r^2s^2}}\,dv\notag \\
&= \frac{1}{4\pi t} e^{-i\frac{r^2+s^2}{4t}} H_0^{(1)}\big(\tfrac{-rs}{2t}\big) \label{eq:2dim}
\end{align}
in the case $n=2$.  Here, $H_0^{(1)}$ is the Hankel function introduced in Section~2.  Some relevant
properties of this function are that $H_0^{(1)}(z)$ is analytic in the upper halfplane and decays
asymptotically like $\sqrt{\pi i/2z}e^{iz}$ as $z \to \infty$ along any ray.

In the case $n=3$, the integral in \eqref{eq:convolution} becomes
\begin{equation*} \begin{aligned}
\int_\R &e^{it\lambda}R_3^-(\lambda,r)R_3^-(\lambda,s)\,d\lambda \\
&= \frac{-2\pi\,e^{-i\frac{r^2+s^2}{4t}}}{(-4\pi i)^3t^2}\int_{2rs}^\infty
{\textstyle \Big[\bigg(\frac{2r^2+v+\sqrt{v^2-4r^2s^2}}
           {r\big(\sqrt{v+2rs}+\sqrt{v-2rs}\big)}\Big)
+ \Big(\frac{2r^2+v-\sqrt{v^2-4r^2s^2}}
             {r\big(\sqrt{v+2rs}-\sqrt{v-2rs}\big)}\Big)\Big]} \\
&\hskip 3.7in\times \frac{e^{-i\frac{v}{4t}}}{\sqrt{v^2-4r^2s^2}}\,dv \\
&= \frac{-2\pi\,e^{-i\frac{r^2+s^2}{4t}}}
        {(-4\pi i)^3\,t^2}\Big(\frac{r+s}{rs}\Big)
  \int_{2rs}^\infty \frac{e^{-i\frac{v}{4t}}}{\sqrt{v-2rs}}\,dv. \\
\end{aligned}
\end{equation*}
At this point it remains to calculate the Fourier transform of an inverse
square-root function, which yields
\begin{equation*}
\int_{2rs}^\infty \frac{e^{-i\frac{v}{4t}}}{\sqrt{v-2rs}}\,dv
  = e^{-i\frac{rs}{2t}} \sqrt{-4\pi i\,t}.
\end{equation*}
The final result is
\begin{equation} \label{eq:3dim}
\int_\R e^{it\lambda}R_3^-(\lambda,r)R_3^-(\lambda,s)\,d\lambda
 = \frac{1}{2i(-4\pi i\,t)^{3/2}}\Big(\frac{r+s}{rs}\Big)
e^{-i\frac{(r+s)^2}{4t}}.
\end{equation}

\subsection{Dimensions $n > 3$}

The recurrence relation for $R_{n+2}^-(\lambda)$ makes it possible to
compute the analogous terms in dimensions $n=5, 7, \ldots$, by
repeatedly applying the differential operator $(4\pi^2 rs)^{-1}
 \frac{\partial^2}{\partial r\partial s}$ to the three-dimensional result
\eqref{eq:3dim}.
For small values of $t$, the leading-order term occurs when all derivatives
fall on $e^{-i(r+s)^2/(4t)}$.  This leads to the following asymptotic
expression as $t\to 0$, which is valid in any odd dimension $n \ge 3$.
\begin{equation} \label{eq:ndim}
\int_\R e^{it\lambda}R_n^-(\lambda,r)R_n^-(\lambda,s) \, d\lambda
 = \frac{1}{2i(-4\pi i\,t)^{n-\frac32}}
  \bigg[\frac{(r+s)^{n-2}}{(rs)^{\frac{n-1}{2}}}\bigg]
e^{-i\frac{(r+s)^2}{4t}}
  + O\big(t^{-(n-\frac52)}\big).
\end{equation}
The same result is true in even dimensions as well.  To see this, recall that $H_0^{(1)}(z) = (\frac{\pi
i}{2z})^{1/2}e^{iz}\omega(z)$, where the derivatives of $\omega$ satisfy the following bounds as $|z|$
goes to infinity:
\begin{equation*}
\lim_{z\to\infty} \omega(z) = 1, \qquad
\big({\textstyle\frac{d}{dz}}\big)^k \omega(z) = O(|z|^{-k}),\ k = 1,2,\ldots
\end{equation*}
The expression in \eqref{eq:2dim} can then be rewritten as
\begin{equation*}
\int_\R e^{it\lambda}R_2^-(\lambda,r)R_2^-(\lambda,s)\,d\lambda
= \frac{1}{2i(-4\pi i\,t\;rs)^{1/2}}e^{-i\frac{(r+s)^2}{4t}}
\omega\big({\textstyle -\frac{rs}{2t}}\big).
\end{equation*}
Applying the differential operators $\ptl[][r]$ and $\ptl[][s]$ only increases
the
degree of the singularity at $t=0$ when the derivative falls on the
term $e^{-i(r+s)^2/(4t)}$.  If the derivative falls instead on
$\omega(-\frac{rs}{2t})$,
one power of $t$ is added to the denominator, but the effect is cancelled by
the faster decay of $\frac{d}{dz}\omega(z)$.
Consequently, when $(4\pi^2 rs)^{-1}\frac{\partial^2}{\partial r \partial s}$
is applied iteratively to \eqref{eq:2dim}, the leading-order term results
from having all of the derivatives fall on $e^{-i(r+s)^2/(4t)}$.
The recurrence relation for $R_{n+2}^-(\lambda)$ then dictates that
\begin{equation} \tag{\ref{eq:ndim}}
\int_\R e^{it\lambda}R_n^-(\lambda,A)R_n^-(\lambda,B)\,d\lambda
 = \frac{1}{2i(-4\pi i\,t)^{n-\frac32}}
  \bigg[\frac{(r+s)^{n-2}}{(rs)^{\frac{n-1}{2}}}\bigg] e^{-i\frac{(r+s)^2}{4t}}
  + O\big(t^{-(n-\frac52)}\big)
\end{equation}
for dimensions $n=4, 6, \ldots$, as desired.  The results of this calculation
can be summarized as follows.

\begin{proposition}
Suppose $n \ge 3$ and let $K$ be a compact subset of $(0,\infty)$.  There exist  constants $C_1, C_2 <
\infty$, depending on $n$ and $K$, such that the remainder function
\begin{equation*} 
G(r,s,t) :=
\int_\R e^{it\lambda}R_n^-(\lambda,A)R_n^-(\lambda,B)\,d\lambda
 - \frac{1}{2i(-4\pi i\,t)^{n-\frac32}}
  \bigg(\frac{(r+s)^{n-2}}{(rs)^{\frac{n-1}{2}}}\bigg) e^{-i\frac{(r+s)^2}{4t}}
\end{equation*}
satisfies the estimates
\begin{equation*}
|G(r,s,t)| \leq C_1 t^{-(n-\frac52)}, \qquad
\Big|\ptl[][t]G(r,s,t)\Big| \leq C_2 t^{-(n-\frac12)},
\end{equation*}
uniformly in $r,s \in K$ and $0<t\leq 1$.
\end{proposition}

\begin{proof}
One obtains an exact expression for $G(r,s,t)$ by differentiating the
base case $n=2$ or
$n=3$.  Under the assumption $r,s \in K$, every monomial in $r$ and $s$
(including those with fractional and/or negative exponents) can be dominated
by a constant.  Every expression of the form
$t^{-k}\frac{d^k}{dz^k}\omega\big(\frac{-rs}{2t}\big)$ can also be bounded
by a constant.  Finally, nonnegative powers of $t$ are smaller than 1.

The function $G(r,s,t)$ consists of all the lower-order terms where at
least one of the partial derivatives $\ptl[][r], \ptl[][s]$ does not fall on
the exponential $e^{-i(r+s)^2/(4t)}$.
It follows that each of these terms is $O(t^{-(n-\frac52)})$.
If the derivative $\ptl[][t]$ is taken at the end, this can only increase
the sharpness of the singularity by a factor of $t^{-2}$.
\end{proof}

To be precise, the proposition above is describing the Fourier transform of a distribution as the
integrand $R_n^-(\lambda,r)R_n^-(\lambda,s)$ experiences growth on the order of
$|\lambda|^{\frac{n-3}2}$.  In Lemma~\ref{lem:asymptotic}, the auxilliary function $\psi_L(\lambda)$ is
introduced to make the integral absolutely convergent.  This has the effect of convolving the
distribution $G(s,t,\cdot)$ with the approximate identity $(2\pi)^{-1}\widehat{\psi_L}$.

At a fixed time $0<t \leq 1$, if $L \geq 2t^{-1}$ one can estimate the
effect of the convolutions
\begin{equation*}
\Big| \big[(2\pi)^{-1}\widehat{\psi_L} * (\cdot)^{-(n-\frac32)}
  e^{-i\frac{r^2+s^2}{4(\cdot)}}\big]\, (t)  - t^{-(n-\frac32)}
   e^{-i\frac{r^2+s^2}{4t}} \Big|
 \leq C_{n,K}L^{-1}t^{-(n+\frac12)}
\end{equation*}
and
\begin{equation*}
\Big| \big[(2\pi)^{-1}\widehat{\psi_L} * G(r,s,\cdot)]\,(t) - G(r,s,t)\Big|
\leq C_{n,K}L^{-1}t^{-(n-\frac12)}
\end{equation*}
by using the Mean Value Theorem and the support
property of $\hat{\psi}$.  If $L > Ct^{-3}$, these resulting differences
are no larger than the initial size estimate for $G(r,s,t)$.  Furthermore,
at fixed $0<t\leq 1$ they vanish in the limit $L\to\infty$ uniformly over
all pairs $r,s\in K$.

Recall the definition of $I_L(t,|x-x_1|,|x_1-y|)$ in the notation of this
section:
\begin{multline*}
I_L(t,|x-x_1|,|x_1-y|) = \int e^{it\lambda} \Big[R_n^+(\lambda,|x-x_1|)
  R_n^+(\lambda,|x_1-y|) \\
- R_n^-(\lambda,|x-x_1|)R_n^-(\lambda,|x_1-y|)\Big]\,d\lambda.
\end{multline*}
Under the substitutions $r = |x-x_1|$ and $s = |x_1-y|$, we have fully characterized the contribution of
the term $e^{it\lambda}R_n^-(\lambda,r)R_n^-(\lambda,s)$ to the integral. The inverse Fourier transform
of $R_n^+(\lambda,r)R_n^+(\lambda,s)$ is a distribution supported on the half line $\{t \leq 0\}$
because of analyticity considerations.  After convolution with $\widehat{\psi_L}$, it will be supported
in $(-\infty,L^{-1}]$ and therefore vanishes at any $t > 0$ once $L > t^{-1}$.

This concludes the proof of Lemma~\ref{lem:asymptotic}.

\end{document}